\documentclass[12pt,a4paper,titlepage,oneside]{book}
\usepackage[english,german]{babel}
\usepackage{amssymb}
\usepackage{amsfonts}
\usepackage{amsthm}
\usepackage{graphicx}
\usepackage{rotating}
\usepackage{lscape}
\usepackage{fancyhdr}
\usepackage{diagrams}
\usepackage{psfrag}

\newcommand{\tstamp}{\today}

\title{\sc On the Topological Tverberg Theorem}

\author{ \\ \\ \sc Diploma Thesis\\ \sc submitted by \\\sc Torsten Sch\"oneborn  \\
\\
\\
\\
\\
\\
Supervised by Prof.~Dr.~G\"unter M.~Ziegler\\
Coreferee Priv.-Doz.~Dr.~Michael Joswig
\\
\\
Institut f\"ur Mathematik, Fakult\"at II, \\
Technische Universit\"at Berlin\\ \\}
\date{Berlin, \tstamp}

\begin{document}
\selectlanguage{english}
\pagenumbering{roman} \maketitle



\newtheorem{thm}{Theorem}[section]
\newtheorem{cor}[thm]{Corollary}
\newtheorem{lem}[thm]{Lemma}
\newtheorem{prop}[thm]{Proposition}
\newtheorem*{TTT}{``Topological Tverberg Theorem''}
\newtheorem*{WindingNumberConjecture}{Winding Number Conjecture}
\newtheorem*{dSkeletonConjecture}{\emph{d}-Skeleton Conjecture}
\newtheorem*{False statement}{False statement}
\newtheorem{conj}[thm]{Conjecture}

\theoremstyle{definition}
\newtheorem{defn}[thm]{Definition}
\newtheorem{ex}[thm]{Example}

\theoremstyle{remark}
\newtheorem{rem}[thm]{Remark}

  \newcommand{\mfd}{ma\-ni\-fold}
  \newcommand{\mfds}{ma\-ni\-folds}
  \newcommand{\mfdwb}{ma\-ni\-fold-with-boundary}
  \newcommand{\mfdswb}{ma\-ni\-folds-with-boundary}
  \newcommand{\twod}{two-di\-men\-sional }
  \newcommand{\td}{three-di\-men\-sional }
  \newcommand{\real}{\mathbb{R}}
  \newcommand{\sphere}{\mathbb{S}}
  \newcommand{\con}{connected}
  \renewcommand{\sc}{simply connected}
  \newcommand{\m}{$\mathcal{M}$}
  \newcommand{\mf}{\mathcal{M}}
  \newcommand{\ttt}{Topological Tverberg Theorem}
  \newcommand{\wnc}{Winding Number Conjecture}
  \newcommand{\dsc}{$d$-Skeleton Conjecture}
  \newcommand{\df}{\Delta^{(d+1)(q-1)}}
  \newcommand{\tpa}{Tverberg partition}
  \newcommand{\tpo}{Tverberg point}
  \newcommand{\pl}{piecewise linear}
  \newcommand{\gp}{general position}
  \newcommand{\wpa}{winding partition}
  \newcommand{\wpo}{winding point}
  \newcommand{\codim}{\mathrm{codim}}
  \newcommand{\w}{winding}
  \newcommand{\defw}{\textbf}
  \newcommand{\dist}{\textrm{dist}}
  \newcommand{\wa}{W_{\neq 0}}
  \newcommand{\ind}{\mathrm{ind}}
  \newcommand{\dy}{$\Delta$-\emph{Y}}
  \newcommand{\yd}{\emph{Y}-$\Delta$}

\selectlanguage{english}
\chapter*{Summary}
Helge Tverberg proved in 1966 that for every \emph{linear} map
from the \linebreak $((d+1)(q-1))$-dimensional simplex $\df$ into
$\real^d$ there is a set of $q$ disjoint faces of this simplex
such that their images intersect in a point \cite{tve66}.

It is conjectured that such a set of disjoint faces exists for
every \emph{continuous} map $\df\to\real^d$ as well, but no
complete proof of such a ``Topological Tverberg Theorem'' is known
yet. Up to now, it has only been proven that the conjecture holds
in the case that $q$ is a prime power \cite{vol96}. A proof of the
Topological Tverberg Theorem for arbitrary $q$ is considered as
one of the biggest challenges in topological combinatorics.

Furthermore, it is still unclear if the $q$ disjoint faces in the
Tverberg Theorem can be uniquely determined by the mapping or if
there are several so called Tverberg partitions. Gerard Sierksma
conjectured that for every linear map from the
$((d+1)(q-1))$-dimensional simplex into $\real^d$ there are at
least $((q-1)!)^d$ different Tverberg partitions. He has found an
example that attains exactly that number of partitions, but a
comprehensive proof is not known yet. A lower bound is proven only
for the case that $q$ is prime \cite{vuz93}, but this bound is far
below the bound conjectured by Sierksma.

In Chapter 1, I briefly describe the origin of the Topological
Tverberg Theorem and summarize the current status of research. In
doing so, all fundamental definitions and theorems used in this
thesis are introduced.\\

A linear map from the simplex is determined by its behavior on the
vertices of the simplex. In order to find the Tverberg partitions
of a map, it is therefore enough to know where the vertices of the
simplex are mapped to.

A continuous map on the other hand is not completely determined by
its behavior on a subskeleton. In this thesis, I prove that for
the Topological Tverberg Theorem it suffices to know the mapping
on the $(d-1)$-skeleton. The restriction to this subskeleton is a
considerable simplification, since the dimension of the simplex is
much larger than $d-1$ if $q$ is big. We introduce a ``\wnc'':

\begin{quote}
For every continuous map from the $(d-1)$-skeleton of the
\linebreak $((d+1)(q-1))$-simplex into $\real^d$, we can choose
$q$ disjoint, at most $d$-dimensional faces of the simplex $\df$
together with a point in $\real^d$, such that for every face
either the point is in the image of the face or the image of the
boundary of the face ``winds around'' the point.
\end{quote}

In Chapter 2, I prove that the \wnc\ and the \ttt\ are equivalent.
This is the main theorem of this thesis. The proof is structured
in two parts:

First we point out, that we can restrict the search for \tpa s to
the $d$-skeleton. Therefore, we introduce the \dsc\ and prove its
equivalence to the \ttt. In a second step, we deal with the
equivalence of the \dsc\ and the \wnc. We prove this for the
higher dimensional cases \linebreak ($d\geq3$) first. Afterwards,
we show that the $d$-dimensional case of the \wnc\ follows from
the $d+1$-dimensional case.

At the end of the second chapter, we see that Sierksma's
conjecture about the number of \tpa s transfers to the \wnc.\\

The \wnc\ is particularly intuitive in the case \linebreak $d=2$,
since it deals with maps from the complete graph $K_{3(q-1)+1}$
into the plane. It claims that in every image of this graph either
$q-1$ triangles wind around one vertex or $q-2$ triangles wind
around the intersection of two edges, where the triangles, edges
and vertices are disjoint.

In Chapter 3, we examine which graphs have this property,
especially which subgraphs of $K_{3(q-1)+1}$. The most interesting
result of this chpater is the following: If $q$ is prime, then the
graph $K_{3(q-1)+1}$ has this property even after deleting a
maximal matching. For the case $q=3$, this even constitutes the
minimal
subgraph of $K_7$ having this property.\\
\\

I would like to thank Prof.~Ziegler for the supervision of this
thesis and his many stimulating suggestions. Furthermore, I would
like to thank Stephan Hell and Arnold Wassmer for the motivating
cooperation and their continuing interest in this thesis.

\selectlanguage{german}
\chapter*{Zusammenfassung}
Helge Tverberg bewies 1966, dass es zu jeder \emph{linearen}
Abbildung des\linebreak $((d+1)(q-1))$-dimensionalen Simplex $\df$
in den $\real^d$ eine Menge von $q$ disjunkten Seiten dieses
Simplex gibt, deren Bilder sich in einem Punkt schneiden
\cite{tve66}.
Man vermutet, dass es auch zu jeder \emph{stetigen} Abbildung
$\df\to\real^d$ eine Menge disjunkter Seiten mit dieser
Eigenschaft
gibt, allerdings ist bis heute noch kein vollst{\"a}ndiger Beweis f{\"u}r
ein solches "`Topologisches Tverberg-Theorem"' gefunden worden.
Bislang konnte lediglich gezeigt werden, dass die Vermutung
zutrifft, falls $q$ eine Primzahlpotenz ist \cite{vol96}. Die
G{\"u}ltigkeit des Topologischen Tverberg-Theorems f{\"u}r beliebige $q$
gilt als eine der gr{\"o}{\ss}ten Herausforderungen der topologischen
Kombinatorik.


Ebenfalls unklar ist bislang, ob die $q$ disjunkten Seiten in
Tverbergs Theorem eindeutig festgelegt sein k{\"o}nnen, oder ob es in
Abh{\"a}ngigkeit von $d$ und $q$ f{\"u}r jede Abbildung mehrere solcher
sogenannten Tverberg-Partitionen gibt. Gerard Sierksma vermutet,
dass es f{\"u}r jede lineare Abbildung des \linebreak
$((d+1)(q-1))$-dimensionalen Simplex in den $\real^d$ mindestens
$((q-1)!)^d$ verschiedene Tverberg-Partitionen gibt. Er hat ein
Beispiel gefunden, das genau diese Anzahl an Partitionen hat; ein
allgemeiner Beweis steht aber noch aus. Nur f{\"u}r den Fall, dass $q$
eine Primzahl ist, kennt man bislang eine untere Schranke
\cite{vuz93}. Diese liegt jedoch weit unter der von Sierksma
vermuteten.

In Kapitel 1 beschreibe ich kurz die Entstehung des Topologischen
Tverberg-Theorems und fasse den aktuellen Stand der Forschung
zusammen. Dabei werden alle grundlegenden Definitionen und
Theoreme aufgef{\"u}hrt, die in dieser Arbeit benutzt werden.\\

Eine lineare Abbildung des Simplex ist bestimmt durch die
Abbildung der Ecken des Simplex. Um die Tverberg-Partitionen einer
Abbildung zu finden, gen{\"u}gt es daher zu wissen, wohin die Ecken
des Simplex abgebildet werden.

Im Gegensatz dazu ist eine stetige Abbildung nicht durch ihr
Verhalten auf einem Teilskelett vollst{\"a}ndig festgelegt. In dieser
Arbeit zeige ich, dass es f{\"u}r das Topologische Tverberg-Theorem
dennoch gen{\"u}gt, die Abbildung auf dem $(d-1)$-Skelett zu kennen.
Durch die Beschr{\"a}nkung auf das $(d-1)$-Skelett ergibt sich eine
Vereinfachung, da die Dimension des Simplex f{\"u}r gro{\ss}e $q$ sehr
viel gr{\"o}{\ss}er als $d-1$ ist. Wir stellen die
"`Windungszahlvermutung"' auf:
\begin{quote}
F{\"u}r jede stetige Abbildung vom $(d-1)$-Skelett des
$((d+1)(q-1))$-Simplex in den $\real^d$ k{\"o}nnen wir $q$ disjunkte,
h{\"o}chstens $d$-dimensionale Seiten des Simplex $\df$ und einen
Punkt des $\real^d$ ausw{\"a}hlen, so dass f{\"u}r jede Seite entweder ihr
Bild den Punkt selbst trifft oder aber das Bild ihres Randes den
Punkt ``uml{\"a}uft''.
\end{quote}

In Kapitel 2 beweise ich, dass die Windungszahlvermutung und das
Topologische Tverberg-Theorem {\"a}quivalent sind. Dies ist das
Hauptresultat dieser Arbeit. Der Beweis gliedert sich in zwei
Teile:

Zuerst machen wir uns klar, dass man sich bei der Suche nach
Tverberg-Partitionen auf das $d$-Skelett beschr{\"a}nken kann. Dazu
stellen wir eine\linebreak $d$-Skelett-Vermutung auf und zeigen
ihre {\"A}quivalenz zum Topologischen Tver\-berg-Theorem. Im zweiten
Schritt behandeln wir die {\"A}quivalenz von $d$-Skelett-Vermutung und
Windungszahl-Vermutung. Diese beweisen wir zuerst f{\"u}r die
h{\"o}herdimensionalen F{\"a}lle ($d\geq 3$). Danach zeigen wir, dass aus
dem $d$-dimensionalen Fall der Windungszahl-Vermutung der
$(d-1)$-dimensionale folgt.

Am Schluss des zweiten Kapitels sehen wir, dass sich Sierksmas
Vermutung {\"u}ber die Anzahl der Tverberg-Partitionen auf die
Windungszahlvermutung {\"u}bertragen l{\"a}sst.\\

Besonders anschaulich ist die Windungszahlvermutung im Fall $d=2$,
da sie sich hier mit Abbildungen des vollst{\"a}ndigen Graphen
$K_{3(q-1)+1}$ in die Ebene besch{\"a}ftigt. Sie behauptet, dass in
jedem Bild dieses Graphen entweder eine Ecke von $q-1$ Dreiecken
umlaufen wird oder aber der Schnittpunkt zweier Kanten von $q-2$
Dreiecken umlaufen wird, wobei die Dreiecke, Kanten und Ecken
paarweise disjunkt sind.
In Kapitel 3 untersuchen wir, welche Graphen (insbesondere welche
Teilgraphen des $K_{3(q-1)+1}$) diese Eigenschaft haben. Das
interessanteste Resultat dieses Kapitels ist das folgende: Der
Graph $K_{3(q-1)+1}$ hat sogar abz{\"u}glich eines maximalen Matchings
diese Eigenschaft, falls $q$ eine Primzahl ist. F{\"u}r den Fall $q=3$
ist damit sogar der minimale Teilgraph von $K_7$ mit
dieser Eigenschaft gefunden.\\
\\

Ich danke Herrn Prof.~Ziegler f{\"u}r das interessante Thema, viele
Anregungen und insbesondere f{\"u}r die immer ge{\"o}ffnete T{\"u}r. Arnold
Wassmer und Stephan Hell bin ich dankbar f{\"u}r die motivierende
Zusammenarbeit und das stete Interesse an dieser Arbeit. Euch
beiden und ganz besonders Henryk Gerlach danke ich f{\"u}r das
Korrekturlesen - ohne Euch h{\"a}tte diese Arbeit (noch) mehr Fehler.
Ein besonderer Dank geht an meine liebe Freundin Anna f{\"u}r das
viele Verst{\"a}ndnis, wenn ich mal wieder in Gedanken war. Meinen
Eltern m{\"o}chte ich danken f{\"u}r die Unterst{\"u}tzung w{\"a}hrend des
gesamten Studiums, dessen Abschluss diese Arbeit darstellt.

\selectlanguage{english}

\tableofcontents \clearpage
\chapter{Introduction}
\pagenumbering{arabic}

In this chapter, we first discuss three different versions of the
classical Tverberg Theorem. Then, we introduce the Topological
Tverberg Theorem (which is really a conjecture) and discuss
bounds for the number of \tpa s.

\section{The Tverberg Theorem}

The historical starting point of the topic of this thesis is the
following theorem from linear geometry.

\begin{thm}[Tverberg Theorem]\label{TverbergOriginal}
Let $d$ and $q$ be positive integers. No matter how
$(d+1)(q-1)+1$ points are chosen in $\real^d$, it is always
possible to partition them into $q$ disjoint sets such that the
convex hulls of these sets intersect, i.e., such that they have a
point in common.
\end{thm}

The first proof was delivered by Helge Tverberg \cite{tve66}.
Today, several different ways of proving it are known. Tverberg
himself offered another proof in \cite{tve81}.

By $\Delta^N$ we denote the $N$-dimensional simplex, by
$\Delta^N_k$ its $k$-skeleton. We will normally not distinguish
between a simplicial complex and its realization, unless it could
cause confusion. We can express the Tverberg Theorem in terms of
a linear map:

\begin{thm}[Tverberg Theorem (Equivalent version
I)]\label{TverbergEquivalent1} For every linear map
$$f:\Delta^{(d+1)(q-1)}\rightarrow \real^d$$ there are $q$
disjoint faces of $\df$ such that their images have a point in
common.
\end{thm}

To see that this is equivalent to the original formulation of
Tverberg's theorem, observe that the convex hull of $n$ points in
$\real^d$ is precisely the image of the linear map
$\Delta^{n-1}\to\real^d$ that maps the $n$ vertices of
$\Delta^{n-1}$ to these $n$ points.

\begin{defn}
Let $k$ be a nonnegative integer and let $f$ be a (not
necessarily linear) map $f:\df_k\rightarrow \real^d$ and $S$ a
set of $q$ disjoint faces $\sigma$ of $\df_k$. We call $S$ a
\defw{\tpa} for the map $f$ if the images of the faces in $S$ have a point
in common, that is, if
$$\bigcap_{\sigma\in S}f(\sigma)\neq\emptyset.$$
Every point in this nonempty intersection is called a
\defw{\tpo}. There might be vertices of $\df_k$ that are not
contained in any face of $S$, although this can happen only in
degenerated cases.
\end{defn}

Using this definition, we can formulate Tverberg's theorem even
simpler:

\begin{thm}[Tverberg Theorem (Equivalent version II)]
For every linear map $$f:\Delta^{(d+1)(q-1)}\rightarrow \real^d$$
there is a \tpa.
\end{thm}

Two questions arise now.
\begin{itemize}
\item Does this theorem hold for a wider class of maps, for
example continuous ones, as well?

 \item How many \tpa s are there at least?
\end{itemize}

We will deal with these questions in the following subsections.

\section{The \ttt}

The following conjecture is a generalization of Tverberg's
theorem to arbitrary continuous maps. It is misleadingly referred
to as the ``\ttt'', although up to now no complete proof of this
conjecture is known.

\begin{conj}[``\ttt'']
For every continuous map $$f:\Delta^{(d+1)(q-1)}\rightarrow
\real^d$$ there is a \tpa.
\end{conj}

\begin{rem}
If we want to restrict the parameters $d$ and $q$ in this
conjecture, we for example talk about ``the case $d=4$ of the
\ttt'' or ``the case $q=6$ of the \ttt''.

The \ttt\ was proven in three cases:
\begin{itemize}
\item The case $d=1$ is equivalent to the mean value theorem for
continuous functions $f:\real\to\real$.

\item The \ttt\ for higher dimensions $d$ was first proven for
prime $q$ by B{\'a}r{\'a}ny, Shlosman and Sz{\H{u}}cs \cite{bss81} using
deleted products. A proof using deleted joins and the
$\mathbb{Z}_q$-index is given in \cite{mat03}.

\item \"Ozaydin proved the more general case of $q$ a prime power
1987 in a still unpublished manuscript. Later, Aleksei Volovikov
gave an alternative proof \cite{vol96}. An elaborate version of
Sarkaria's proof using characteristic classes \cite{sak00} can be
found in de Longueville \cite{deL01}.
\end{itemize}
\end{rem}

All other cases still remain open; the smallest open case is
therefore $d=2$, $q=6$ (see also Table
\ref{ProvenCasesOfTheTTT}). This case deals with maps from the
15-dimensional simplex to $\real^2$. In Matou{\v{s}}ek's opinion, ``the
validity of the \ttt\ for arbitrary (nonprime) $q$ is one of the
most challenging problems in this field [topological
combinatorics]'' \cite[p.154]{mat03}. It is known that lower
dimensional cases follow from higher dimensional ones:

\begin{prop}[de Longueville \cite{deL01}]\label{tttFromDToD-1}
If the \ttt\ holds for $q$ and $d$, then it also holds for $q$
and $d-1$.
\end{prop}

\begin{table}
\begin{center}
\begin{tabular}{|c||c|c|c|c|c|c|c|c|c|c|c|c|c|c|c|}
  \hline
  $d\backslash q$ & 1 & 2 & 3 & 4 & 5 & 6 & 7 & 8 & 9 & 10 & 11 & 12 & 13 & 14 & 15\\
  \hline
  \hline
  1 & $\checkmark$ & $\checkmark$ & $\checkmark$ & $\checkmark$ & $\checkmark$ & $\checkmark$ & $\checkmark$ & $\checkmark$ & $\checkmark$ & $\checkmark$ & $\checkmark$ & $\checkmark$ & $\checkmark$ & $\checkmark$ & $\checkmark$\\
  \hline
  2 & $\checkmark$ & $\checkmark$ & $\checkmark$ & $\checkmark$ & $\checkmark$ &   & $\checkmark$ & $\checkmark$ & $\checkmark$ &   & $\checkmark$ &   & $\checkmark$ & & \\
  \hline
  3 & $\checkmark$ & $\checkmark$ & $\checkmark$ & $\checkmark$ & $\checkmark$ &   & $\checkmark$ & $\checkmark$ & $\checkmark$ &   & $\checkmark$ &   & $\checkmark$ & & \\
  \hline
  4 & $\checkmark$ & $\checkmark$ & $\checkmark$ & $\checkmark$ & $\checkmark$ &   & $\checkmark$ & $\checkmark$ & $\checkmark$ &   & $\checkmark$ &   & $\checkmark$ & & \\
  \hline
  5 & $\checkmark$ & $\checkmark$ & $\checkmark$ & $\checkmark$ & $\checkmark$ &   & $\checkmark$ & $\checkmark$ & $\checkmark$ &   & $\checkmark$ &   & $\checkmark$ & & \\
  \hline
  6 & $\checkmark$ & $\checkmark$ & $\checkmark$ & $\checkmark$ & $\checkmark$ &   & $\checkmark$ & $\checkmark$ & $\checkmark$ &   & $\checkmark$ &   & $\checkmark$ & & \\
  \hline
  7 & $\checkmark$ & $\checkmark$ & $\checkmark$ & $\checkmark$ & $\checkmark$ &   & $\checkmark$ & $\checkmark$ & $\checkmark$ &   & $\checkmark$ &   & $\checkmark$ & & \\
  \hline
  8 & $\checkmark$ & $\checkmark$ & $\checkmark$ & $\checkmark$ & $\checkmark$ &   & $\checkmark$ & $\checkmark$ & $\checkmark$ &   & $\checkmark$ &   & $\checkmark$ & & \\
  \hline
  9 & $\checkmark$ & $\checkmark$ & $\checkmark$ & $\checkmark$ & $\checkmark$ &   & $\checkmark$ & $\checkmark$ & $\checkmark$ &   & $\checkmark$ &   & $\checkmark$ & & \\
  \hline
  \vdots & \vdots & \vdots & \vdots & \vdots & \vdots & & \vdots & \vdots & \vdots & & \vdots & & \vdots & & \\

\end{tabular}
\end{center}
\caption{The checkmarks indicate the proven cases of the \ttt
.}
\label{ProvenCasesOfTheTTT}
\end{table}

I will now give an outline of the proof of the \ttt\ in the case
that $q$ is a prime.

\begin{thm}[{B{\'a}r{\'a}ny, Shlosman and Sz{\H{u}}cs \cite{bss81}}]\label{TTTforPrimeQ}
The \ttt\ is valid if $q$ is a prime.
\end{thm}

We follow the proof presented in \cite{mat03}. The central definition is the deleted join,
which we will use as a configuration space. We need two
definitions -- the first one for simplicial complexes, the second
one for topological spaces.

\begin{defn}
Let $q$ be a positive integer. For $q$ sets $A_1,\dots,A_q$, let
$A_1\uplus A_2\uplus \dots\uplus A_q$ be the set
$$(A_1\times\{1\})\cup(A_2\times\{2\})\cup\dots\cup(A_q\times\{q\}).$$

Let $\Delta_1,\Delta_2,\dots,\Delta_q$ be a simplicial complexes
with vertex sets $V_1,V_2,\dots,V_q$. The
\defw{join} $\Delta_1\ast\Delta_2\ast\dots\ast\Delta_q$ has vertex set
$$V_1\uplus V_2\uplus\dots\uplus V_q$$ and
face set $$\{F_1\uplus F_2\uplus \dots\uplus F_q \mid
F_i\mathrm{\ is\ a\ face\ of\ }\Delta_i\mathrm{\ for\ all\
}i\}.$$

Let $t_1,\dots,t_q$ be nonnegative real numbers with
$\sum_{i=1}^qt_i=1$. For all $i$, let $p_i$ be a point in a
realization of $F_i$. We define the following notation for points
in the realization of the face $F_1\uplus F_2\uplus \dots\uplus
F_q$:
$$(t_1p_1\oplus t_2p_2\oplus\dots\oplus t_qp_q):=\sum_{v\in
F_1}t_1t_{1,v}(v\times\{1\})+\dots+\sum_{v\in
F_q}t_qt_{q,v}(v\times\{q\})$$ where $t_{i,v}\geq 0$ such that $p_i=\sum_{v\in
F_i}t_{i,v} v$ and $\sum_{v\in F_i}t_{i,v}=1$ for all $i$.

Let $\Delta$ be a simplicial complex with vertex set $V$. The
\defw{\emph{q}-fold join} $\Delta^{*q}$ of $\Delta$ is
$$\Delta^{*q}:=\underbrace{\Delta\ast\Delta\ast\dots\ast\Delta}_{q\mathrm{\ times}}.$$

The \defw{\emph{q}-fold pairwise deleted join}
$\Delta^{*q}_{\Delta(2)}$ is the subcomplex of $\Delta^{*q}$ with
the face set
$$\{F_1\uplus F_2\uplus \dots\uplus F_q \mid \mathrm{the\ }F_i
\mathrm{\ are\ pairwise\ disjoint\ faces\ of\ }\Delta \}.$$
\end{defn}

\begin{defn}
Let $X_1,X_2,\dots,X_q$ be topological spaces. The
\defw{join} \linebreak $X_1\ast X_2\ast\dots\ast X_q$ is the topological
space
$$X_1\ast X_2\ast\dots\ast X_q:=X_1\times X_2\times \dots\times X_q\times Y/\approx$$
where $Y$ is the convex hull of the standard unit vectors in
$\real^q$, that is, the set
$$\{(t_1,\dots,t_q)\in\real^q\mid \sum_{i=1}^q t_i=1 \mathrm{\
and\ } t_i\geq 0 \mathrm{\ for\ all\ }i\},$$ and $\approx$ is
given by
\begin{eqnarray*}
& (p_1,\dots,p_q,(t_1,\dots,t_q))\approx
(p'_1,\dots,p'_q,(t'_1,\dots,t'_q))\\
& \mathrm{if\ and\ only\ if\ }(t_i=t'_i)\mathrm{\ and\ }(t_i\neq
0\Rightarrow p_i=p'_i)\mathrm{\ for\ all\ }i.
\end{eqnarray*}

We write
$$(t_1p_1\oplus t_2p_2\oplus\dots\oplus t_qp_q):=(p_1,\dots,p_q,(t_1,\dots,t_q)).$$

Let $X$ be a topological space. The \defw{\emph{q}-fold join}
$X^{*q}$ of $X$ is the topological space
$$\underbrace{X\ast X\ast\dots\ast X}_{q\mathrm{\
times}}.$$

The \defw{\emph{q}-fold \emph{q}-wise deleted join}
$X^{*q}_{\Delta}$ is the following subspace of $X^{*q}$:
\begin{eqnarray*}
X^{*q}_{\Delta}& := &
X^{*q}\backslash\{(p,\dots,p,(\textstyle\frac{1}{q},\dots,\frac{1}{q}))\mid p\in X\}\\
& = &
X^{*q}\backslash\{(\textstyle\frac{1}{q}p\oplus\dots\oplus\frac{1}{q}p)\mid
p\in X\}
\end{eqnarray*}
\end{defn}

\begin{defn}
Let $f:\Delta\to X$ be a continuous map from the realization of a
simplicial complex $\Delta$ to a topological space $X$. We define
the \defw{\emph{q}-fold join} $f^{*q}$ of $f$ to be the following
map.
\begin{eqnarray*}
f^{*q} : \Delta^{*q} & \to & X^{*q}\\
f^{*q}(t_1p_1\oplus t_2p_2\oplus\dots\oplus t_qp_q) & := &
(t_1f(p_1)\oplus t_2f(p_2)\oplus\dots\oplus t_qf(p_q))
\end{eqnarray*}
\end{defn}

The definitions of deleted joins are tailored for our purposes:
Let us assume that a counterexample $f:\df\to\real^d$ exists.
First, we take its $q$-fold join
$$f^{*q}:(\df)^{*q}\to(\real^d)^{*q}$$
and restrict it to the deleted join
$$f^{*q}|_{(\df)^{*q}_{\Delta(2)}}:(\df)^{*q}_{\Delta(2)}\to(\real^d)^{*q}.$$
Observe that every face of $(\df)^{*q}_{\Delta(2)}$ represents
$q$ (possibly empty) disjoint faces of $\df$. Because we assumed
$f$ to be a counterexample, these faces do not form a \tpa. Thus
no point $(t_1p\oplus\dots\oplus t_qp)$ of $(\real^d)^{*q}$ is in
the image $f^{*q}|_{(\df)^{*q}_{\Delta(2)}}$, where $p$ is in
$\real^d$ and the $t_i$ are positive real numbers. In particular
no point $(\frac{1}{q}p\oplus\dots\oplus\frac{1}{q}p)$ is
attained. Therefore we can reduce the co-domain to
$(\real^d)^{*q}_\Delta$ and obtain a map
$$f^{*q}_\Delta:(\df)^{*q}_{\Delta(2)}\to(\real^d)^{*q}_\Delta.$$

If we show that there is no such map $f^{*q}_\Delta$, then we can
immediately deny the existence of a counterexample for the \ttt.
Of course, a map $(\df)^{*q}_{\Delta(2)}\to(\real^d)^{*q}_\Delta$
exists (for example the constant map), but $f^{*q}_\Delta$ has an
important property: It is $\mathbb{Z}_q$-equivariant.
$\mathbb{Z}_q$ operates on the $q$-fold join $A^{*q}$ of a
simplicial complex or topological space $A$ by
$$\mu.(t_1p_1\oplus t_2p_2\oplus\dots\oplus t_qp_q):=(t_2p_2\oplus t_3p_3\oplus\dots\oplus t_qp_q\oplus
t_1p_1)$$ where $\mu$ is a generator of $\mathbb{Z}_q$. Note that
$\mathbb{Z}_q$ operates in the same way on deleted joins.
Furthermore, the $q$-fold join of a map is
$\mathbb{Z}_q$-equivariant under this operation. As we have seen,
the question whether an arbitrary map
$(\df)^{*q}_{\Delta(2)}\to(\real^d)^{*q}_\Delta$ exists is not
fruitful; but does a \emph{$\mathbb{Z}_q$-equivariant} map
$(\df)^{*q}_{\Delta(2)}\to(\real^d)^{*q}_\Delta$ exist?

The answer is ``no'', if $q$ is a prime. We prove this using
index theory.

\begin{defn}
Let $X$ be a topological space or a simplicial complex with a
$\mathbb{Z}_q$-operation. We define the
\defw{$\mathbb{Z}_q$-index} of $X$ with respect to this operation
as
$$\ind_{\mathbb{Z}_q}(X):=\min\{n\mid\mathrm{there\ is\ a\ }\mathbb{Z}_q\mathrm{-equivariant\ map\
}X\to(\Delta^n)^{*q}_{\Delta(2)}\}$$ Here, we regard
$(\Delta^n)^{*q}_{\Delta(2)}$ equipped with the
$\mathbb{Z}_q$-operation described above.
\end{defn}

\begin{lem}\label{IndexTheoryLemma1}
Let $X$ and $Y$ be spaces with $\mathbb{Z}_q$-operation. There is
no \linebreak $\mathbb{Z}_q$-invariant map $X\to Y$ if
$\ind_{\mathbb{Z}_q}(X)$ is greater than $\ind_{\mathbb{Z}_q}(Y)$.
\end{lem}

\begin{proof}
Assume an equivariant map $f:X\to Y$ exists. By definition, there
is also an equivariant map $g:Y\to
(\Delta^{\ind_{\mathbb{Z}_q}(Y)})^{*q}_{\Delta(2)}$. By combining
these, we obtain an equivariant map
$$g\circ f:X\to (\Delta^{\ind_{\mathbb{Z}_q}(Y)})^{*q}_{\Delta(2)}.$$
Again by definition we conclude that
$\ind_{\mathbb{Z}_q}(X)\leq\ind_{\mathbb{Z}_q}(Y)$, which
contradicts the conditions of the lemma.
\end{proof}

\begin{lem}[see Matou{\v{s}}ek {\cite[Sections 6.2 and 6.3]{mat03}}]\label{IndexTheoryLemma2} \-
\begin{itemize}
\item Let $X$ and $Y$ be spaces with $\mathbb{Z}_q$-operation.
$\mathbb{Z}_q$ operates on $X*Y$ by operating on both factors
simultaneously  and we have
$$\ind_{\mathbb{Z}_q}(X*Y)\leq\ind_{\mathbb{Z}_q}(X)+\ind_{\mathbb{Z}_q}(Y)+1.$$

\item Let $\mathbb{Z}_q$ operate on $\mathbb{S}_1$ by a rotation
of angle $\frac{2\pi}{q}$. Then we have
$$\ind_{\mathbb{Z}_q}(\mathbb{S}^1)=1.$$

\item The index of a pairwise deleted join of a simplex is
$$\ind_{\mathbb{Z}_q}((\Delta^n)^{*q}_{\Delta(2)})=n.$$

\item If $q$ is a prime, then the index of
$((\real^d)^{*q}_\Delta)$ is
$$\ind_{\mathbb{Z}_q}((\real^d)^{*q}_\Delta)=(d+1)(q-1)-1.$$
\end{itemize}
\end{lem}

\begin{proof}[Proof of Theorem \ref{TTTforPrimeQ}]
If a counterexample existed, then there would be an equivariant
map
$f^{*q}_\Delta:(\df)^{*q}_{\Delta(2)}\to(\real^d)^{*q}_\Delta$.
But there is no such equivariant map, since
$$\ind_{\mathbb{Z}_q}((\df)^{*q}_{\Delta(2)})=(d+1)(q-1)>(d+1)(q-1)-1=\ind_{\mathbb{Z}_q}((\real^d)^{*q}_\Delta).$$
\end{proof}

\section{How many \tpa s are there?}

Sierksma conjectured that for every linear map $f:\df\to\real^d$
there are at least $((q-1)!)^d$ \tpa s. This number is attained
for the configuration of $d+1$ tight clusters, with $q-1$ points
each, placed at the vertices of a simplex, and one point in the
middle.

For $d=1$, the mean value theorem implies Sierksma's conjecture.
In almost all other cases, Sierksma's conjecture is still
unresolved at the time of writing (see Table
\ref{SierksmaNumberTable}). Nevertheless, for special values of
$q$, a lower bound is known:

\begin{thm}[Vu{\v{c}}i{\'c} and {\v{Z}}ivaljevi{\'c} \cite{vuz93}]\label{ManyTverbergPartitions}
If $q$ is a prime, then there are at least
$$\frac{1}{(q-1)!}\cdot\Big(\frac{q}{2}\Big)^{(d+1)(q-1)/2}$$
\tpa s for every continuous map $f:\df\to\real^d$.
\end{thm}

A nice proof can also be found in Matou{\v{s}}ek {\cite[Theorem
6.6.1]{mat03}}. For arbitrary $q$ but linear $f:\df\to\real^d$,
the best known lower bound is 1, given by the classical Tverberg
Theorem. Furthermore, no non-topological method is known to yield
a good lower bound.

\begin{landscape}
\begin{table}
\begin{tabular}{|c||rr|rr|rr|rr|rr|rr|}
  \hline
$d\backslash q$ & \multicolumn{2}{c|}{1} & \multicolumn{2}{c|}{2} & \multicolumn{2}{c|}{3} & \multicolumn{2}{c|}{4} & \multicolumn{2}{c|}{5} & \multicolumn{2}{c|}{6} \\
\hline\hline
1 &   1 & (1) & 1 & (1) & 2     & (2)  & 6         &  (6) &  24          & (24)        &  120         & (120) \\
2 &   1 & (1) & 1 & (1) & 4     & (2)  & 36        &  (1) &  576         & (11)        &  14400       & (0)   \\
3 &   1 & (1) & 1 & (1) & 8     & (3)  & 216       &  (1) &  13824       & (64)        &  1728000     & (0)   \\
4 &   1 & (1) & 1 & (1) & 16    & (4)  & 1296      &  (1) &  331776      & (398)       &  2,07$\cdot 10^{8}$     & (0)   \\
5 &   1 & (1) & 1 & (1) & 32    & (6)  & 7776      &  (1) &  7962624     & (2484)      &  2,48$\cdot 10^{10}$    & (0)   \\
6 &   1 & (1) & 1 & (1) & 64    & (9)  & 46656     &  (1) &  1,91$\cdot 10^{8}$     & (15523)     &  2,98$\cdot 10^{12}$    & (0)   \\
7 &   1 & (1) & 1 & (1) & 128   & (13) & 279936    &  (1) &  4,58$\cdot 10^{9}$     & (97013)     &  3,58$\cdot 10^{14}$    & (0)   \\
8 &   1 & (1) & 1 & (1) & 256   & (20) & 1679616   &  (1) &  1,10$\cdot 10^{11}$    & (606330)    &  4,29$\cdot 10^{16}$    & (0)   \\
9 &   1 & (1) & 1 & (1) & 512   & (29) & 10077696  &  (1) &  2,64$\cdot 10^{12}$    & (3789562)   &  5,15$\cdot 10^{18}$    & (0)   \\
10&   1 & (1) & 1 & (1) & 1024  & (44) & 60466176  &  (1) &  6,34$\cdot 10^{13}$    & (23684758)  &  6,19$\cdot 10^{20}$    & (0)   \\
11&   1 & (1) & 1 & (1) & 2048  & (65) & 362797056 &  (1) &  1,52$\cdot 10^{15}$    & (148029737) &  7,43$\cdot 10^{22}$    & (0)   \\
\hline \multicolumn{5}{c}{}\\ \hline
$d\backslash q$ & \multicolumn{4}{c|}{7} & \multicolumn{2}{c|}{8} & \multicolumn{2}{c|}{9} & \multicolumn{2}{c|}{10} & \multicolumn{2}{c|}{11} \\
\hline\hline
1 &   \multicolumn{2}{r}{720}         & \multicolumn{2}{r|}{(720)}         & 5040       & (5040) &  40320     & (40320) &  362880   & (362880)  & 3628800  & (3628800) \\
2 &   \multicolumn{2}{r}{518400}      & \multicolumn{2}{r|}{(110)}         & 25401600   & (1)    &  1,62$\cdot 10^{9}$   & (1)     &  1,31$\cdot 10^{11}$ & (0)       & 1,31$\cdot 10^{13}$ & (35130)   \\
3 &   \multicolumn{2}{r}{3,73$\cdot 10^{8}$}     & \multicolumn{2}{r|}{(4694)}        & 1,28$\cdot 10^{11}$   & (1)    &  6,55$\cdot 10^{13}$  & (1)     &  4,77$\cdot 10^{16}$ & (0)       & 4,77$\cdot 10^{19}$ & (1,76$\cdot 10^{8}$) \\
4 &   \multicolumn{2}{r}{2,68$\cdot 10^{11}$}    & \multicolumn{2}{r|}{(201228)}      & 6,45$\cdot 10^{14}$   & (1)    &  2,64$\cdot 10^{18}$  & (1)     &  1,73$\cdot 10^{22}$ & (0)       & 1,73$\cdot 10^{26}$ & (8,89$\cdot 10^{11}$)\\
5 &   \multicolumn{2}{r}{1,93$\cdot 10^{14}$}    & \multicolumn{2}{r|}{(8627646)}     & 3,25$\cdot 10^{18}$   & (1)    &  1,06$\cdot 10^{23}$  & (1)     &  6,29$\cdot 10^{27}$ & (0)       & 6,29$\cdot 10^{32}$ & (4,47$\cdot 10^{15}$)\\
6 &   \multicolumn{2}{r}{1,39$\cdot 10^{17}$}    & \multicolumn{2}{r|}{(3,69$\cdot 10^{8}$)}     & 1,63$\cdot 10^{22}$   & (1)    &  4,29$\cdot 10^{27}$  & (1)     &  2,28$\cdot 10^{33}$ & (0)       & 2,28$\cdot 10^{39}$ & (2,25$\cdot 10^{19}$)\\
7 &   \multicolumn{2}{r}{1,00$\cdot 10^{20}$}    & \multicolumn{2}{r|}{(1,58$\cdot 10^{10}$)}    & 8,26$\cdot 10^{25}$   & (1)    &  1,73$\cdot 10^{32}$  & (1)     &  8,28$\cdot 10^{38}$ & (0)       & 8,28$\cdot 10^{45}$ & (1,13$\cdot 10^{23}$)\\
8 &   \multicolumn{2}{r}{7,22$\cdot 10^{22}$}    & \multicolumn{2}{r|}{(6,79$\cdot 10^{11}$)}    & 4,16$\cdot 10^{29}$   & (1)    &  6,98$\cdot 10^{36}$  & (1)     &  3,00$\cdot 10^{44}$ & (0)       & 3,00$\cdot 10^{52}$ & (5,70$\cdot 10^{26}$)\\
9 &   \multicolumn{2}{r}{5,19$\cdot 10^{25}$}    & \multicolumn{2}{r|}{(2,91$\cdot 10^{13}$)}    & 2,09$\cdot 10^{33}$   & (1)    &  2,81$\cdot 10^{41}$  & (1)     &  1,09$\cdot 10^{50}$ & (0)       & 1,09$\cdot 10^{59}$ & (2,87$\cdot 10^{30}$)\\
10&   \multicolumn{2}{r}{3,74$\cdot 10^{28}$}    & \multicolumn{2}{r|}{(1,25$\cdot 10^{15}$)}    & 1,05$\cdot 10^{37}$   & (1)    &  1,13$\cdot 10^{46}$  & (1)     &  3,95$\cdot 10^{55}$ & (0)       & 3,95$\cdot 10^{65}$ & (1,44$\cdot 10^{34}$)\\
11&   \multicolumn{2}{r}{2,69$\cdot 10^{31}$}    & \multicolumn{2}{r|}{(5,35$\cdot 10^{16}$)}    & 5,33$\cdot 10^{40}$   & (1)    &  4,57$\cdot 10^{50}$  & (1)     &  1,43$\cdot 10^{61}$ & (0)       & 1,43$\cdot 10^{72}$ & (7,27$\cdot 10^{37}$)\\
\hline
\end{tabular}
\caption{The number of \tpa s conjectured by Sierksma. The number
in brackets shows the currently highest proven lower bound for
continuous maps.} \label{SierksmaNumberTable}
\end{table}
\end{landscape}

\chapter{The \wnc}

We saw in the introduction that there are two equivalent versions
of the Tverberg Theorem for linear maps: Theorem
\ref{TverbergEquivalent1}, that deals with maps $f$ from the
entire $((d+1)(q-1))$-dimensional Simplex to $\real^d$, and the
original version (Theorem \ref{TverbergOriginal}), that uses only
the image of the vertices of the simplex, that is, it talks about
$f|_{\df_0}$. In this chapter, we will establish conjectures that
are equivalent to the \ttt, but talk only about $f|_{\df_k}$ for
some integer $k$ less than $(d+1)(q-1)$.

In fact, we will see two such conjectures: The \dsc\linebreak
($k=d$) claims that every continuous map $f:\df_d\to\real^d$ has
a \tpa\ in the $d$-skeleton. We will derive the equivalence to
the \ttt\ by reducing the problem to maps in ``\gp''.

The \wnc\ ($k=d-1$) claims that for every continuous map
$f:\df_{d-1}\to\real^d$, the boundary of many simplices of
$\df_d$ wind about a point contained in the image of some of the
faces. We also consider bounds on the number of \tpa s
respectively \wpa s.

We start our discussion with the \wnc, since it is more powerful,
and introduce the \dsc\ later as a link between the \wnc\ and the
\ttt. We need the following definition.

\begin{defn}
First, let us assume $d\geq 2$. Let
$f:\mathbb{S}^{d-1}\rightarrow\real^d$ be a continuous map from
the $(d-1)$-dimensional sphere to $\real^d$, and let $p$ be a
point in $\real^d$. We choose an isomorphism
$I:\pi_{d-1}(\real^d\backslash\{p\})\rightarrow\mathbb{Z}$. If
$f$ does not attain $p$, that is, if $p\notin
f(\mathbb{S}^{d-1})$, then we define the
\defw{winding number of \emph{f} with respect to \emph{p}} as
$$W(f,p):=I([f])\in\mathbb{Z}.$$
The sign of $W(f,p)$ depends on the choice of $I$, but the
expression $$``W(f,x)=0"$$ is independent of this choice. Let
$\partial\Delta^d:=\Delta^d_{d-1}$. For maps
$f:\partial\Delta^d\to\real^d$, we define $W(f,p)$ via
$h:\mathbb{S}^{d-1}\stackrel{\approx}{\to}\Delta^d_{d-1}=\partial\Delta^d$
to be $W(f,x):=W(f\circ h,x)$.

Now if $d=1$, then we do not have such an isomorphism $I$. In
this thesis, we are not interested in the exact value of
$W(f,p)$, but only in whether it is zero. Therefore it is
sufficient for our purposes to define $W(f,p)$ to be zero if the
two points $f(\mathbb{S}^0)$ lie in the same component of
$\real\backslash \{p\}$. Otherwise we say that $W(f,p)\neq 0$.
\end{defn}

\begin{WindingNumberConjecture}\label{wnc}
For all positive integers $d$ and $k$ and every continuous map
$f:\Delta^{(d+1)(q-1)}_{d-1}\rightarrow \real^d$ there are $q$
disjoint faces $\sigma_1,\dots,\sigma_q$ of
$\Delta^{(d+1)(q-1)}_d$ and a point $p\in\real^d$ such that for
every $i$, one of the following holds:
 \begin{itemize}
 \item $\dim(\sigma_i)\leq d-1$ and $p\in f(\sigma_i)$
 \item $\dim(\sigma_i)= d$ and ($p\in f(\partial\sigma_i)$ or $W(f|_{\partial\sigma_i},p)\neq
 0$)
 \end{itemize}
Such a set $S=\{\sigma_1,\dots,\sigma_q\}$ will be called a
\defw{\wpa}; $p$ will be called a \defw{\wpo}.
\end{WindingNumberConjecture}

We intentionally included the case ``$p\in f(\partial\sigma_i)$''
in the previous conjecture. See also Remark
\ref{WhyTheCondition}.

\begin{figure*}
\begin{center}
\includegraphics*[bb= 110 550 470 670,scale=.7]{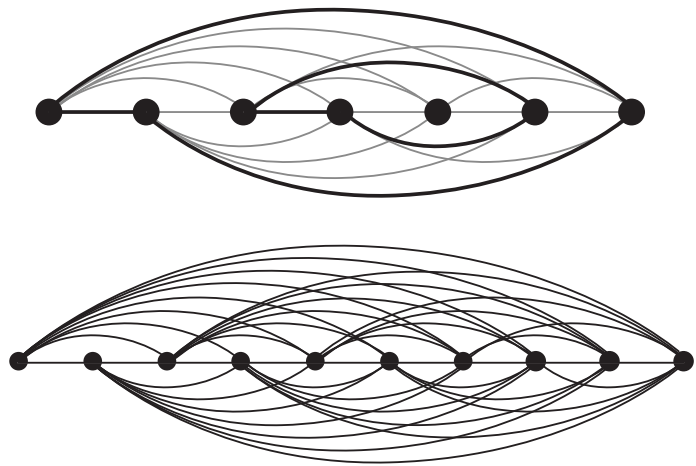}
\includegraphics*[bb= 100 390 500 530,scale=.7]{K_nAlsLinieMitZweiKanten.eps}
\end{center}
\caption{The alternating linear model of $K_7$ and $K_{10}$. The
thick black lines in the drawing of $K_7$ form a \wpa.}
\label{K_n}
\end{figure*}

\begin{ex}\label{K_nByGuy}
Let us look at a concrete example of a continuous map
$\df_{d-1}\to\real^d$. In the case $d=2$, this is really a
drawing of $K_{3(q-1)+1}$, the complete graph with
$3(q-1)+1=3q-2$ vertices. If the drawing is in ``\gp'' (in a way
made precise in the next section), then the \wnc\ says that in
the drawing of $K_{3q-2}$ either $q-1$ (possibly distorted)
triangles wind around one vertex, or $q-2$ triangles wind around
the intersection of two edges, with the triangles, edges and the
vertex being pairwise disjoint.

The following way to draw $K_n$ is called ``the alternating
linear model'' and was proposed in \cite{saa64}. Draw the $n$
vertices on a line and number them from left to right. Draw the
edges $[i,i+1]$ on this line and the edges $[i,i+k],k\geq 2$ on
one side of the line (e.g.~above) if $i$ is odd and on the other
side (e.g.~below) if $i$ is even. Figure \ref{K_n} illustrates
the situation for $K_7$ and $K_{10}$.

For this drawing the \wnc\ is satisfied: The vertex with number
$2q-1$ is a \wpo. For example the $q-1$ disjoint triangles
$\langle 1,2,3q-2\rangle,\langle 3,4,3q-3\rangle,\dots,\langle
2q-3,2q-2,3q-q\rangle$ wind around it.

This is not surprising, because the classical Tverberg Theorem
guarantees that every rectilinear drawing of $K_{3q-2}$ satisfies
the \wnc, and there is a rectilinear drawing for the alternating
linear model.
\end{ex}

We can give an alternative description of the term ``\wpa'':
\begin{defn}
For any map $f:\partial\Delta^d\to\real^d$, we define
$$\wa(f):=f(\partial\Delta^d)\cup\{x\in\real^d\backslash f(\partial\Delta^d)\mid W(f,x)\neq 0\}.$$

\end{defn}

\begin{rem}\label{WhyTheCondition}
Later on, it will be advantageous that $\wa(f)$ is a closed set
containing $f(\partial\Delta^d)$ even in degenerated cases where
$\wa(f)$ might be empty. This is why we had to add
$f(\partial\Delta^d)$ to the definition of $\wa(f)$ and include
``$p\in f(\partial\sigma_i)$'' in our first formulation of the
\wnc.
\end{rem}

\begin{lem}
A set $S=\{\sigma_1,\dots,\sigma_q\}$ of $q$ disjoint faces of
$\df_d$ is a \wpa\ for $f:\Delta^{(d+1)(q-1)}_{d-1}\rightarrow
\real^d$ if and only if
$$\bigcap_{\dim(\sigma_i)<d} f(\sigma_i)
\cap\bigcap_{\dim(\sigma_i)=d}\wa(f|_{\partial\sigma_i})\neq\emptyset.$$
\begin{flushright}
\vskip -.7cm $\Box$
\end{flushright}
\end{lem}

\begin{WindingNumberConjecture}[Equivalent version]
For every continuous map $f:\Delta^{(d+1)(q-1)}_{d-1}\rightarrow
\real^d$ there are $q$ disjoint faces $\sigma_1,\dots,\sigma_q$
of $\Delta^{(d+1)(q-1)}_d$ such that
$$\bigcap_{\dim(\sigma_i)<d} f(\sigma_i)
\cap\bigcap_{\dim(\sigma_i)=d}\wa(f|_{\partial\sigma_i})\neq\emptyset.$$
\end{WindingNumberConjecture}

This conjecture can be proved easily if $d=1$ (see Proposition
\ref{wncFuerD=1}). The rest of this chapter covers the proof of
the following theorem.

\begin{thm}
The \wnc\ is equivalent to the \ttt.
\end{thm}

\begin{figure}
\includegraphics*[bb= 5 370 590 800,scale=.7]{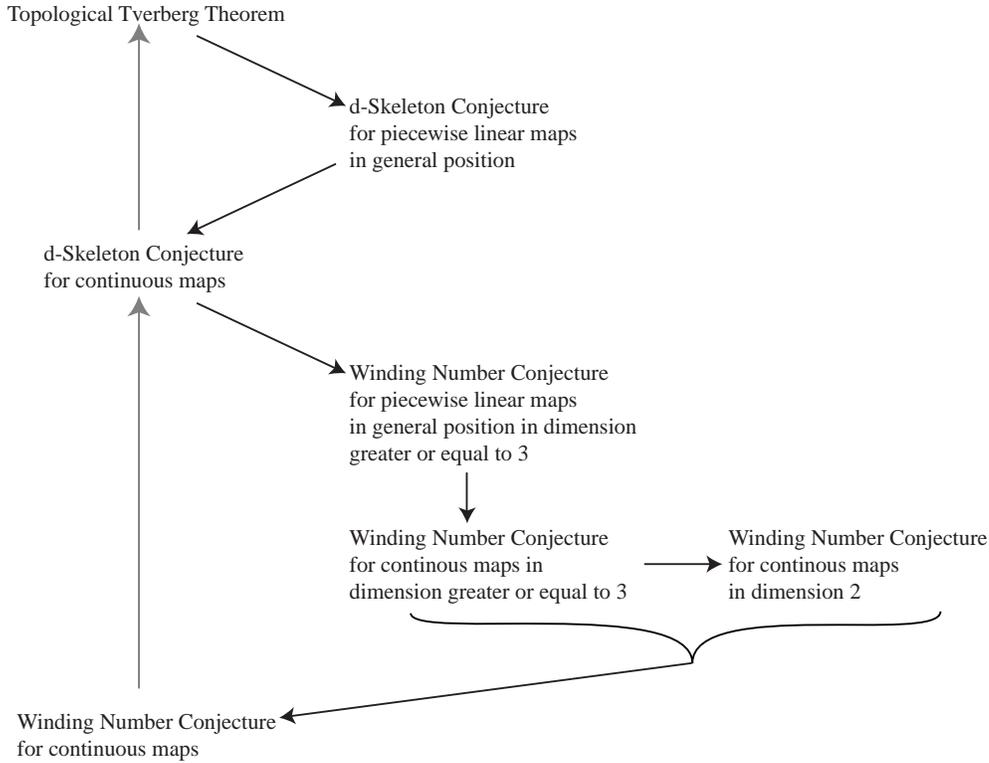}
\caption{A flow chart of the implications between conjectures
proved in this chapter. The grey arrows are the obvious
implications.} \label{flowchart}
\end{figure}

The line of argument of the proof is illustrated in Figure
\ref{flowchart}.

\begin{rem}\label{BasicIdea}
The basic idea of the proof are the following two speculations.
\begin{itemize}
\item Let $F:\df\to\real^d$ be a continuous map. Every \wpa\ of
$F|_{\df_{d-1}}$ is a \tpa\ of $F$.

\item Let $f:\df_{d-1}\to\real^d$ be a continuous map. Then $f$
can be extended to a continuous map $F:\df\to\real^d$ such that
every \tpa\ of $F$ is a \wpa\ of $f$.
\end{itemize}
The first statement turns out to be true, but the second one
needs some adjustment, as we will see in the course of the proof.
We will come back to the above speculations in Theorem
\ref{ConnectionDSCandWNC}.
\end{rem}

\section{Step 1: Reduction to the \emph{d}-skeleton}
First, we show that the \ttt\ guarantees the existence of a \tpa\
in the $d$-skeleton of $\df$.

\begin{dSkeletonConjecture}
Every continuous map $f:\df_d\rightarrow \real^d$ has a \tpa.
\end{dSkeletonConjecture}

\begin{prop}\label{ttttodsc}
The \dsc\ is equivalent to the \ttt.
\end{prop}

It is obvious that the \dsc\ implies the \ttt. The converse is
harder. Its proof is the aim of this subsection. We divide the
proof into the Lemmas \ref{dscpart1} and \ref{dscpart2}.

\subsection{Maps in \gp}

For the first lemma, we need the following definition.

\begin{defn}
Let $\Delta$ be a simplicial complex. A map
$f:\Delta\rightarrow\real^d$ is \defw{linear} if it is linear on
every face of $\Delta$. Such an linear map $f$ is in
\defw{\gp} if for every set of disjoint faces
$\{\sigma_1,\sigma_2,\dots,\sigma_k\}$ of $\Delta$ the inequality
$$\codim(\bigcap^{k}_{i=1}f(\sigma_i)) \geq
\sum^k_{i=1}\codim(f(\sigma_i))$$ holds, where
$\codim(\sigma):=d-\dim(\sigma)$ if $\sigma\subset\real^d$. We
use the convention $\dim(\emptyset)=-\infty$ and thus set
$\codim(\emptyset)=\infty$. The last equation includes the case
$\bigcap_{i=1}^k f(\sigma_i)=\emptyset$. Thus, in that case the
condition above holds independently of the right hand side
because we have $\codim(\bigcap_{i=1}^k f(\sigma_i))=\infty$.
\end{defn}


\begin{figure}
\begin{center}
\includegraphics*[bb= 50 550 285 785,scale=.5]{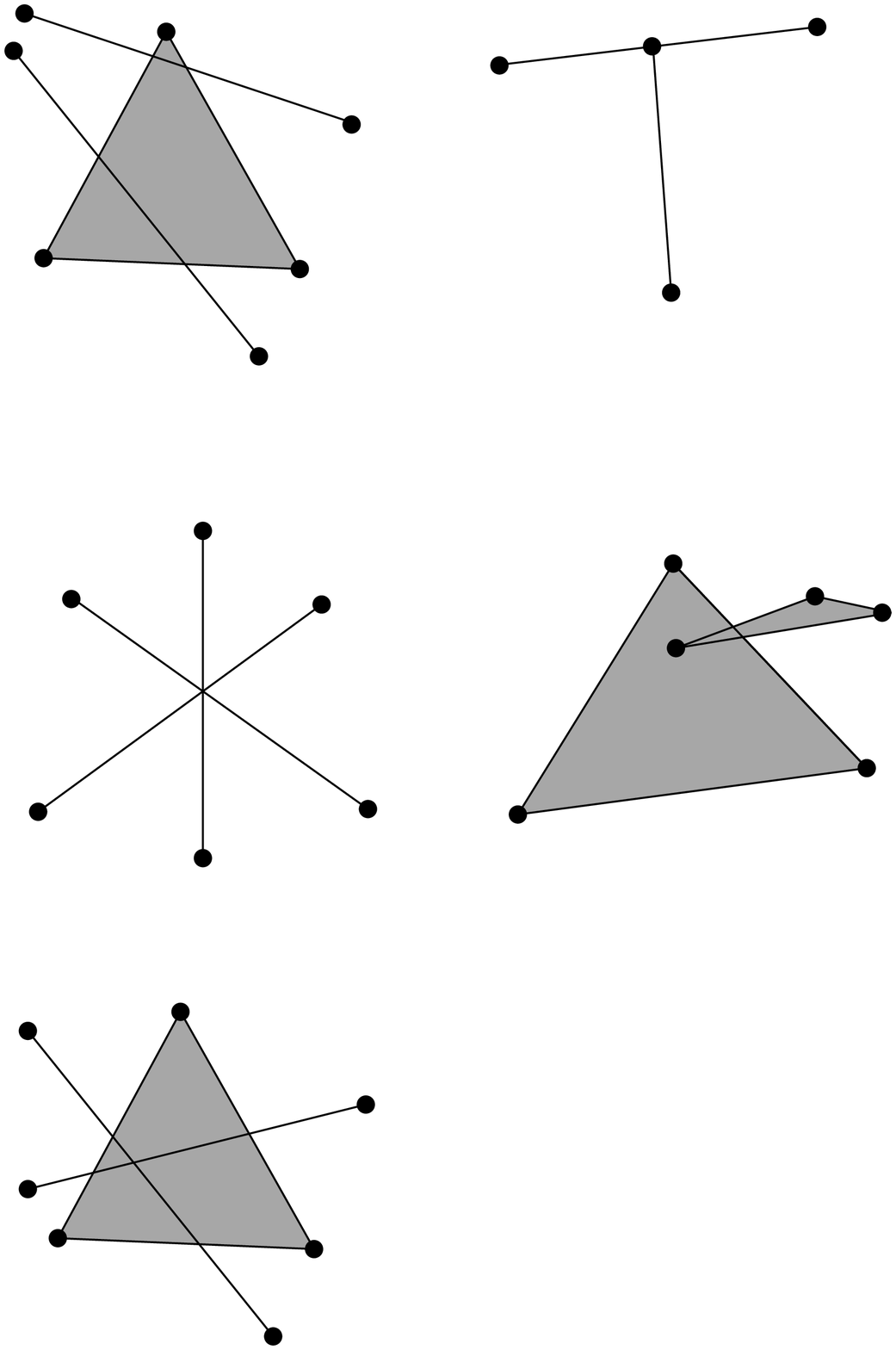}
\includegraphics*[bb= 330 270 570 520, scale=.5]{GeneralPositionExamplesAndCounterexamples.eps}
\includegraphics*[bb= 50 30 300 250, scale=.5]{GeneralPositionExamplesAndCounterexamples.eps}
\end{center}
\caption{Images $f(\Delta)$ of linear maps $f:\Delta\to\real^2$
in \gp.}
\end{figure}

\begin{figure}
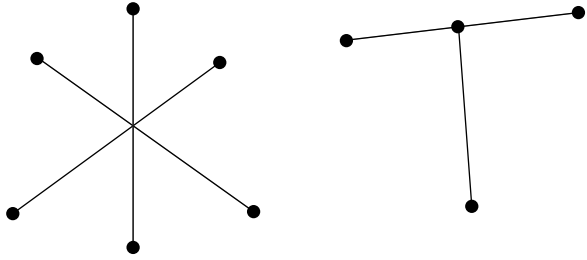

\begin{center}
\includegraphics*[bb= 40 270 300 520, scale=.5]{GeneralPositionExamplesAndCounterexamples.eps}
\includegraphics*[bb= 310 550 560 785,scale=.5]{GeneralPositionExamplesAndCounterexamples.eps}
\end{center}
\caption{Images $f(\Delta)$ of linear maps $f:\Delta\to\real^2$
\emph{not} in \gp. In the last picture, the complex $\Delta$
consists of two lines.}
\end{figure}

We need to restrict ourselves to \pl\ maps to exclude ``wild''
maps.

\begin{defn}
Let $\Delta$ be a simplicial complex. A map
$f:\Delta\rightarrow\real^d$ is \defw{\pl} if there is a
subdivision $s:\Delta'\to\Delta$ such that the composition
$f\circ s:\Delta'\to\real^d$ is a linear map. Furthermore, we
call $f$
\defw{in \gp} if we can choose the subdivision $s$ such that the
linear map $f\circ s$ is in \gp.
\end{defn}

\begin{figure}
\begin{center}
\includegraphics*[bb= 40 620 570 730, scale=.65]{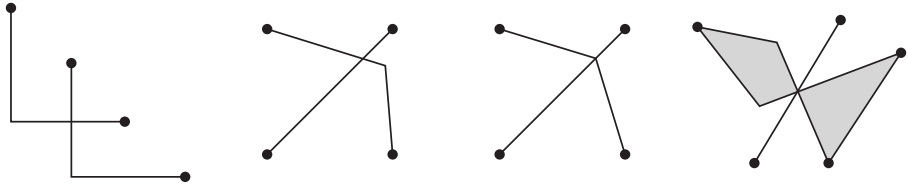}
\end{center}
\caption{Images $f(\Delta)$ of \pl\ maps $f:\Delta\to\real^2$. In
the first three pictures, $\Delta$ consists of two lines, in the
last picture $\Delta$ consists of a triangle and a line. The two
pictures on the left are in \gp, the two on the right are
not.}\label{PLexamples}
\end{figure}

Whether $f\circ s$ is in \gp\ depends on the subdivision $s$. For
example, the  map $f$ depicted on the very left in Figure
\ref{PLexamples} combined with the second barycentric subdivision
gives a linear map not in \gp, although $f$ itself is in \gp.

The definition of \gp\ made here may seem overly restrictive for
the purpose of this section, but we need it in Proposition
\ref{wnckgeq3}.

The key point of maps in \gp\ is the following lemma.

\begin{lem}
Let $\Delta$ be a simplicial complex and
$f:\Delta\rightarrow\real^d$ a \pl\ map in \gp. If
$\{\sigma_1,\sigma_2,\dots,\sigma_k\}$ is a set of disjoint faces
of $\Delta$, then we have
$$\codim(\bigcap^{k}_{i=1}f(\sigma_i)) \geq
\sum^k_{i=1}\max(0,(d-\dim(\sigma_i))).$$
\end{lem}

\begin{rem}
$\bigcap_{i=1}^k f(\sigma_i)$ might have parts of different
dimension. We use the convention $$\dim(A\cup
B):=\max(\dim(A),\dim(B))$$ and thus $$\codim(A\cup
B):=\min(\codim(A),\codim(B))$$
\end{rem}

\begin{proof}
Let $s:\Delta'\to\Delta$ be a subdivision such that $f\circ s$ is
a linear map in \gp.
\begin{eqnarray*}
\codim(\bigcap^{q}_{i=1}f(\sigma_i)) & = &
\min_{\stackrel{\tilde{\sigma}_i\subset\sigma_i}{\tilde{\sigma}_i\mathrm{\ Simplex\ in\ }\Delta'}}\codim(\bigcap^{q}_{i=1}f\circ s(\tilde{\sigma}_i))\\
& \geq & \min_{\tilde{\sigma}_i\subset\sigma_i}\sum^q_{i=1}\codim(f\circ s(\tilde{\sigma}_i))\\
& = & \min_{\tilde{\sigma}_i\subset\sigma_i}\sum^q_{i=1}(d-\dim(f\circ s(\tilde{\sigma}_i)))\\
& = & \sum^q_{i=1}\min_{\tilde{\sigma}_i\subset\sigma_i}(d-\dim(f\circ s(\tilde{\sigma}_i)))\\
& = & \sum^q_{i=1}(d-\max_{\tilde{\sigma}_i\subset\sigma_i}\dim(f\circ s(\tilde{\sigma}_i)))\\
& \geq & \sum^q_{i=1}(d-\min(d,\dim\sigma_i)) \\
& = & \sum^k_{i=1}\max(0,(d-\dim(\sigma_i))).
\end{eqnarray*}
\end{proof}

We need an approximation lemma to tackle continuous maps.

\begin{lem}[Piecewise Linear Approximation Lemma]
Let $\Delta$ be a simplicial complex with a subcomplex
$\Delta_0\subset \Delta$ and let $\varepsilon>0$. Furthermore,
let $f:\Delta\to\real^d$ be a continuous map that is \pl\ on
$\Delta_0\subset \Delta$. Then there is a \pl\ map
$\tilde{f}:\Delta\to\real^d$ that equals $f$ on $\Delta_0$ and
approximates it on the rest of $\Delta$, that is,
$\tilde{f}|_{\Delta_0}=f|_{\Delta_0}$ and
$$\|\tilde{f}-f\|_{\infty}=\max\{|\tilde{f}(x)-f(x)|:x\in\Delta\}<\varepsilon.$$
\end{lem}

\begin{proof}
Let $s:\Delta_0'\to\Delta_0$ be a subdivision such that
$f|_{\Delta_0}\circ s$ is a linear map. We can extend $s$ to a
subdivision $S:\Delta'\to\Delta$ of all of $\Delta$. Since
$\Delta$ is compact, there is an iterated barycentric subdivision
$\tilde{S}:\tilde{\Delta}\to\Delta'$ of $\Delta'$ with the
following property: If $p\in\tilde{\Delta}$ is contained in the
face $\{v_1,\dots,v_k\}$ of $\tilde{\Delta}$, then
$\|f(p)-f(v_i)\|<\varepsilon$. Let
$\tilde{f}:\tilde{\Delta}\to\real^d$ be the linear map that is
given on the vertices $v$ of $\tilde{\Delta}$ by
$\tilde{f}(v):=f(S(\tilde{S}(v)))$. Therefore $\tilde{f}$ is
piecewise linear on $\Delta$, equals $f$ on $\Delta_0$ (because
already $f\circ S$ is linear on this subcomplex) and approximates
$f$ on the rest of $\Delta$.
\end{proof}

\subsection{\tpa s in the \emph{d}-skeleton}

Using the Approximation Lemma and the properties of \pl\ maps, we
can now prove that if the \ttt\ holds, then the \dsc\ holds as
well.

\begin{lem}\label{dscpart1}
Every \tpa\ of any \pl\ map \linebreak $f:\df\to\real^d$ in \gp\
contains only faces of dimension at most $d$.
\end{lem}

\begin{cor}\label{dscpart1b}
If the \ttt\ is true, then the \linebreak \dsc\ holds for all
\pl\ maps in general position.
\end{cor}

\begin{lem}\label{dscpart2}
For every continuous map $f:\df_d\rightarrow\real^d$ there is an
$\varepsilon_f>0$ such that the following holds: If
$\tilde{f}:\df_d\rightarrow\real^d$ is a continuous map with
$\|\tilde{f}-f\|_{\infty}<\varepsilon_f,$ then every \tpa\ of
$\tilde{f}$ is also a \tpa\ of $f$.
\end{lem}

This lemma states that by distorting $f$ by less than
$\varepsilon_f$, we do not create new \tpa s.

\begin{cor}\label{dscpart3}
If the \dsc\ holds for all \pl\ maps in general position, then it
is true in general (i.e., for all continuous maps).
\end{cor}

By these two corollaries, the \ttt\ implies the \dsc.

\begin{proof}[Proof of Lemma \ref{dscpart1}]
Let $f$ be in \gp\ that has an arbitrary \tpa\
$\{\sigma_1,\sigma_2,\dots,\sigma_q\}$.
\begin{eqnarray*}
d & \stackrel{(1)}{\geq} & \codim(\bigcap^{q}_{i=1}f(\sigma_i)) \\
& \stackrel{(2)}{\geq} & \sum^q_{i=1}\max(0,(d-\dim\sigma_i))\\
& \stackrel{(\ast)}{\geq} & \sum^q_{i=1}(d-\dim\sigma_i)\\
& = & qd-(\sum^q_{i=1}((\textrm{number of vertices of
}\sigma_i)-1))\\
& \geq & qd-(\textrm{(number of vertices of }\df)-q)\\
& = & qd-((d+1)(q-1)+1-q)\\
& = & d.
\end{eqnarray*}

\renewcommand{\labelenumi}{(\theenumi):}
\begin{enumerate}
\item This holds because $\{\sigma_1,\sigma_2,\dots,\sigma_q\}$
is a \tpa\ and thus $\bigcap^{q}_{i=1}f(\sigma_i)\neq\emptyset$.

\item This holds because $f$ is in \gp.
\end{enumerate}
In $(\ast)$, equality holds only if $d-\dim(\sigma_i)\geq 0$, or
equivalently if $\dim(\sigma_i)\leq d$ for all $i$, which is what
we had to prove.
\end{proof}

\begin{proof}[Proof of Corollary \ref{dscpart1b}]
Let us assume we are given a \pl\ map $f:\df_d\rightarrow
\real^d$ in \gp. First, we extend this map piecewise linearly to
a map $\tilde{F}:\df\rightarrow\real^d$. This is possible by
taking a continuous extension, which exists because $\real^d$ is
contractible, and then taking a \pl\ approximation without
changing the map on $\df_d$. Now we can obtain a \pl\ map
$F:\df\to\real^d$ in \gp\ by perturbing $\tilde{F}$, again
leaving it unchanged on $\df_d$.

By the \ttt, there is a \tpa\ for~$F$, which must lie in the
$d$-skeleton by Lemma \ref{dscpart1} and thus be a \tpa\ of $f$
as well.
\end{proof}

\begin{proof}[Proof of Lemma \ref{dscpart2}]
Let $\tilde{f}:\df_d\rightarrow \real^d$ be a map that satisfies
$\|\tilde{f}-f\|_{\infty}<\varepsilon$. We want to show that if
$\varepsilon$ is sufficiently small, then we can be sure that
every \tpa\ of $\tilde{f}$ is a \tpa\ of $f$.

If $S$ is a \tpa\ of $\tilde{f}$, then we obtain
$$\emptyset\ \ \neq\ \ \bigcap_{\sigma\in S}\tilde{f}(\sigma)\ \ \subseteq\ \
\bigcap_{\sigma\in S}\{x\in\real^d\mid
\dist(x,f(\sigma))\leq\varepsilon\}.$$

The last expression denotes a compact set that gets smaller when
$\varepsilon$ decreases. If it is empty for $\varepsilon=0$, then
it must therefore already be empty for a sufficiently small
$\varepsilon$, that is, there is a $\varepsilon_S>0$ such that
$$\bigcap_{\sigma\in S}f(\sigma)=\emptyset\ \ \Rightarrow\ \
\bigcap_{\sigma\in S}\{x\in\real^d\mid
\dist(x,f(\sigma))\leq\varepsilon_S\}=\emptyset.$$ If we choose
$$\varepsilon:=\varepsilon_f:=\min\{\varepsilon_S\mid S\textrm{ a set of } q\textrm{ disjoint faces of
}\df\},$$ then we can be sure that every \tpa\ of $\tilde{f}$ is
also a \tpa\ of $f$.
\end{proof}

\begin{proof}[Proof of Corollary \ref{dscpart3}]
For all $\varepsilon > 0$ and continuous maps
$f:\df_d\rightarrow\real^d$, we can obtain a \pl\ map $\tilde{f}$
in general position by distorting $f$ by less than $\varepsilon$.
This distortion can, for example, be carried out via a \pl\
approximation to obtain a \pl\ map followed by a small adjustment
to get this map into \gp.

We restrict the distortion and adjustment to $\varepsilon_f$,
that is, we make sure that $\|\tilde{f}-f\|<\varepsilon_f$. By
assumption, $\tilde{f}$ has a \tpa. This is also a \tpa\ of $f$
(Lemma \ref{dscpart2}).
\end{proof}

\subsection{The connection between \tpa s in the full simplex and in the \emph{d}-skeleton}

In the previous section, we proved the equivalence of the \ttt\
and the \dsc. The arguments that we saw establish the following
stronger result as well.

\begin{prop}\label{ConnectionTTTandDSC}
Let $F:\df\to\real^d$ be a continuous map. Every \tpa\ of
$F|_{\df_d}$ is also a \tpa\ of $F$.

Let $f:\df_d\to\real^d$ be a continuous map. We can extend a
slightly distorted version of $f$ to a continuous map
$F:\df\to\real^d$ such that every \tpa\ of $F$ is also a \tpa\ of
$f$.
\end{prop}

\begin{proof}
The first part is obvious. For the second part, first approximate
$f$ by a \pl\ map $\bar{f}:\df_d\to\real^d$ in \gp\ that is
sufficiently close to $f$ (Lemma \ref{dscpart2}) and extend
$\bar{f}$ to $F$ in the way described in the proof of Corollary
\ref{dscpart1b}.
\end{proof}

\begin{cor}
The \dsc\ is valid if $d=1$ and if $q$ is a prime power.
\end{cor}

Therefore Table \ref{ProvenCasesOfTheTTT} also applies to the
\dsc.

\section{Step 2: Reduction to the (\emph{d}--1)-skeleton}
Now we proceed to prove the equivalence of the \wnc\ and the \dsc.

\begin{prop}\label{WNCtoDSC}
The \wnc\ implies the \dsc.
\end{prop}

\begin{proof}
Let $f:\df_d\to\real^d$ be a continuous map and
$\sigma_1,\dots,\sigma_q$ be a \wpa\ for $f|_{\df_{d-1}}$ with
\wpo\ $P\in\real^d$. This \wpa\ is also a \tpa\ for $f$:
\begin{itemize}
\item If $\dim(\sigma_i)\leq d-1$, we have $P\in
f|_{\df_{d-1}}(\sigma_i)=f(\sigma_i)$.

\item If $\dim(\sigma_i)=d$, then $W(f|_{\partial\sigma_i},P)\neq
0$, hence $P\in f(\sigma_i)$.
\end{itemize}
\end{proof}

The proof of the converse is harder. For this we want to show
that any map $$f:\df_{d-1}\rightarrow\real^d$$ can be extended to
a map
$$F:\df_d\rightarrow\real^d$$ such that every \tpa\ of $F$ is
a \wpa\ of $f$. This would be easy to do if for each
$d$-dimensional face $\sigma\subset\df$, we could satisfy
$F(\sigma)\subset\wa(f|{\partial\sigma})$. Unfortunately, this is
not always possible because of the following proposition:

\begin{prop}
Not every continuous map $f:\mathbb{S}^{d-1}\rightarrow\real^d$
is nullhomotopic within $\wa(f)$.
\end{prop}

We look at two examples. Let $\mathbb{B}^d$ be the
$d$-dimensional ball with boundary $\sphere^{d-1}$.

\begin{figure}
\begin{center}
\includegraphics[220,420][390,580]{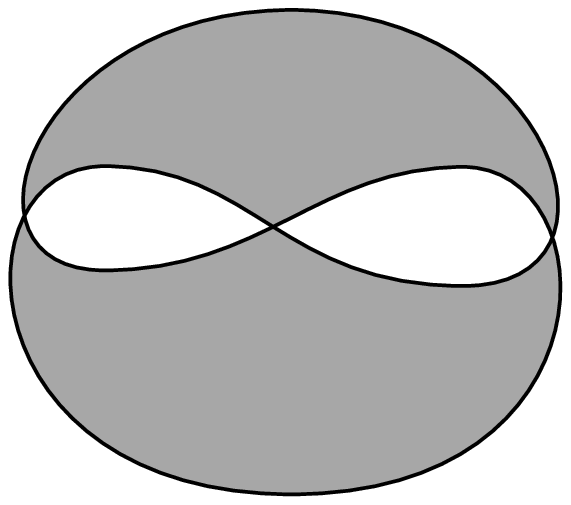}
\end{center}
\caption{A map $f:\sphere^1\to\real^2$ that is not nullhomotopic
within $\wa(f)$. The shaded area is $\wa(f)$.}
\label{d2counterexample}
\end{figure}

\begin{ex}
The first counterexample is a map $f:\sphere^1\to\real^2$
illustrated by Figure \ref{d2counterexample}. The topological
space $\wa(f)$ is homotopy equivalent to the wedge of two spheres
$\sphere^1$. The fundamental group $\pi_1(\wa(f))$ is therefore
equal to
$\pi_1(\sphere^1\vee\sphere^1)=\mathbb{Z}\ast\mathbb{Z}$, the
free product of $\mathbb{Z}$ with itself. The element
$[f]\in\pi_1(\wa(f))$ can be written as the nonzero term
$aba^{-1}b^{-1}$ if we choose generators $a,b$ of
$\mathbb{Z}\ast\mathbb{Z}$ suitably.

If we extend $f$ to $\mathbb{B}^2$, then its image covers at
least one of the two ``holes'' in $\wa(f)$ entirely, which is a
$2$-dimensional set. There is no one-dimensional subset
$V\subset\real^2$ such that $f$ is contractible in $\wa(f)\cup V$.

The suspension $\sphere f:\sphere^2\to\real^3$ is not a
counterexample. We have $$\wa(\sphere
f)=\sphere\wa(f)=\sphere(\sphere^1\vee\sphere^1)=\sphere^2\vee\sphere^2$$
again, but this time the homotopy group
$\pi_2(\sphere^2\vee\sphere^2)$ is not a free product but a free
sum $\mathbb{Z}\oplus\mathbb{Z}$, therefore we calculate
$[f]=aba^{-1}b^{-1}=aa^{-1}bb^{-1}=0$ in $\pi_2(\wa(\sphere f))$.
\end{ex}

\begin{ex}
For $d\geq 4$, the homotopy group $\pi_{d-1}(\sphere^{d-2})$ is
nontrivial (see Hatcher \cite[Sections 4.1 and 4.2]{hat02}). For
example, the Hopf map $\sphere^3\to\sphere^2$ is not
nullhomotopic. Choose such a map
$f:\sphere^{d-1}\to\sphere^{d-2}$ that is not nullhomotopic. Let
$i:\sphere^{d-2}\to\real^d$ be an embedding of the sphere in a
$(d-1)$-dimensional linear subspace of $\real^d$. Then
$\wa(i\circ f)=i(f(\sphere^{d-2}))$, hence $i\circ f$ can not be
contracted in $\wa(i\circ f)$.

An important difference between this example and the previous one
is that here, $i \circ f$ can be contracted within the
$(d-1)$-dimensional subspace that contains $i(\sphere^{d-2})$. No
$d$-dimensional set outside of $\wa(i\circ f)$ is necessary for a
contraction.
\end{ex}

Because of these problems, we have to take a more technical route.

\subsection{Reduction to \pl\ maps in \gp}

We need an approximation lemma similar to \ref{dscpart2}.

\begin{lem}\label{wncpart2}
For every continuous map $f:\df_{d-1}\rightarrow\real^d$ there is
an $\varepsilon_f>0$ such that the following holds: If
$\tilde{f}:\df_{d-1}\rightarrow\real^d$ is a continuous map with
$\|\tilde{f}-f\|_{\infty}<\varepsilon_f$, then every \wpa\ of
$\tilde{f}$ is also a \wpa\ of $f$.
\end{lem}

\begin{cor}\label{wncPiecewiseSufficient}
If the \wnc\ holds for \pl\ maps in \gp, then it also holds for
all continuous maps.
\end{cor}

\begin{proof}[Proof of Lemma \ref{wncpart2}]
Let $\tilde{f}:\df_{d-1}\rightarrow \real^d$ be a map that
satisfies $\|\tilde{f}-f\|_{\infty}<\varepsilon$. We want to show
that if $\varepsilon$ is sufficiently small, then we can be sure
that every \wpa\ of $\tilde{f}$ is a \wpa\ of $f$.

If $S$ is a \wpa\ of $\tilde{f}$, we obtain

\begin{eqnarray*}
\emptyset & \neq & \Big(\bigcap_{\sigma\in
S,\dim(\sigma)<d}\tilde{f}(\sigma)\Big)\cap\Big(\bigcap_{\sigma\in
S,\dim(\sigma)=d}\wa(\tilde{f}|_{\partial\sigma})\Big) \\
& \subseteq & \Big(\bigcap_{\sigma\in
S,\dim(\sigma)<d}\{x\in\real^d\mid
\dist(x,f(\sigma))\leq\varepsilon\}\Big)\\
& & \cap\ \ \Big(\bigcap_{\sigma\in
S,\dim(\sigma)=d}\{x\in\real^d\mid
\dist(x,\wa(f|_{\partial\sigma}))\leq \varepsilon\}\Big).
\end{eqnarray*}

The last expression denotes a compact set that gets smaller when
$\varepsilon$ decreases. If it is empty for $\varepsilon=0$, then
it is
empty for a sufficiently small $\varepsilon_S$.
If we choose
$$\varepsilon_f:=\min\{\varepsilon_S\mid S\textrm{ a set of disjoint faces of
}\df_d\},$$ then we can be sure that every \wpa\ of $\tilde{f}$
is also a \wpa\ of $f$.
\end{proof}

\begin{proof}[Proof of Corollary \ref{wncPiecewiseSufficient}]
Identical to the proof of Corollary \ref{dscpart3}.
\end{proof}

\subsection{The case \emph{d} $\geq$ 3}

\begin{defn}
A \defw{triangulation of $\real^d$} is a simplicial complex
$\Delta$ with a fixed linear map
$\|\Delta\|\stackrel{\cong}{\rightarrow}\real^d$. We do not
distinguish between a face of the triangulation and the
corresponding set in $\real^d$.

Let $\Delta_1,\Delta_2,\dots,\Delta_\ell$ be triangulations of
$\real^d$. They are \defw{in \gp\ with respect to each other} if
for every subset $S\subset \{1,\dots,\ell\}$ and faces
$\sigma_{i}$ of $\Delta_i$, we have
$$\codim(\bigcap_{i\in S}\sigma_i) \geq \sum_{i\in
S}\codim(\sigma_i).$$
\end{defn}

\begin{prop}\label{wnckgeq3}
Let $k\geq3$. If the \dsc\ is true for $d=k$, then the \wnc\
holds for $d=k$.
\end{prop}

\begin{proof}
By Corollary \ref{wncPiecewiseSufficient}, we can restrict
ourselves to \pl\ maps. Let $f:\df_{d-1}\rightarrow\real^d$ be a
\pl\ map in \gp. We divide the proof in three steps:
\begin{enumerate}
\item Choose a triangulation $\Delta_\sigma$ of $\real^d$ for
every face $\sigma\subset\df_d$.

\item Extend $f:\df_{d-1}\to\real^d$ to a continuous map
$F:\df_d\to\real^d$ ``compatible'' to the $\Delta_\sigma$.

\item Show that every \tpa\ of $F$ is a \wpa\ of $f$. By the
\dsc, $F$ has a \tpa, that thus is a \wpa\ for $f$.
\end{enumerate}

\emph{Step 1:} For every face $\sigma\subset\df_d$ choose a
triangulation $\Delta_{\sigma}$ of $\real^d$ such that
\begin{itemize}
\item for $\dim(\sigma)\leq d-1$, the set $f(\sigma)$ is a subset
of the $\dim(\sigma)$-skeleton of $\Delta_\sigma$ and

\item for $\dim(\sigma)=d$, the set $f(\partial\sigma)$ is a
subset of the $(d-1)$-skeleton of $\Delta_\sigma$.
\end{itemize}

Choose the $\Delta_{\sigma}$ such that if
$\sigma_1,\dots,\sigma_\ell$ are disjoint faces of $\df_d$, then
$\Delta_{\sigma_1},\dots,\Delta_{\sigma_\ell}$ are in \gp\ with
respect to each other. This is possible because $f$ is in \gp.
(Here we need the restrictive definition of ``\gp''!)
Furthermore, choose them such that for every $\Delta_\sigma$
there is a map
$b_\sigma:\mathbb{S}^{d-1}\rightarrow\Delta_{\sigma}(=\real^d)$
with the following properties:
\begin{itemize}
\item The map $b_{\sigma}$ is a simplicial embedding with respect
to a suitably chosen triangulation of $\mathbb{S}^{d-1}$.

\item The set $f(\sigma)$ respectively $f(\partial\sigma)$ is
completely contained in the bounded component $B_\sigma$ of
$\real^d\backslash b_\sigma(\mathbb{S}^{d-1})$.

\end{itemize}

\emph{Step 2:} Now, we extend $f$ to a $d$-face
$\sigma\subset\df_d$. Let $\sigma_1,\dots,\sigma_k$ be the
$d$-faces of $B_\sigma \cap \overline{(\real^d\backslash
\wa(f|_{\partial\sigma}))}$ with respect to the triangulation
$\Delta_\sigma$. For every $i$ in $\{1,\dots,k\}$, choose a point
$x_i$ in $\sigma_i$. Let $$\imath_i:B_\sigma\backslash
\{x_1,\dots,x_k\}\hookrightarrow\real^d\backslash\{x_i\}$$ be the
inclusion and let
$I_i:\pi_{d-1}(\real^d\backslash\{x_i\})\stackrel{\cong}{\to}\mathbb{Z}$
be the isomorphism used for the definition of the winding number
$W(\cdot,x_i)$. We have the following commutative diagram:
\begin{diagram}
\pi_{d-1}(B_\sigma\backslash(\stackrel{\circ}{\sigma}_1\cup\dots\cup\stackrel{\circ}{\sigma}_k)) & & & \bigoplus_{i=1}^k\mathbb{Z}\\
\uTo_{(r_1)_\ast}^{\stackrel{(1)}{\cong}} & & & \uTo_{(I_1,\dots,I_k)}^{\stackrel{(5)}{\cong}}\\
\pi_{d-1}(B_\sigma\backslash\{x_1,\dots,x_k\}) & \rTo^{((\imath_1)_\ast,\dots,(\imath_k)_\ast)} & & \bigoplus_{i=1}^k \pi_{d-1}(\real^d\backslash\{x_i\})\\
\dTo_{(r_2\circ\jmath^{-1})_\ast}^{\stackrel{(2)}{\cong}} & & & \uTo_{(\jmath_1)_\ast\oplus\dots\oplus(\jmath_k)_\ast}^{\stackrel{(4)}{\cong}}\\
\pi_{d-1}(\bigvee_{i=1}^k\sphere^{d-1}) & \lTo^{\tilde{\jmath}}_{\stackrel{(3)}{\cong}} & & \bigoplus_{i=1}^k\pi_{d-1}(\sphere^{d-1})
\end{diagram}

\renewcommand{\labelenumi}{(\theenumi):}
\begin{enumerate} \item By blowing up the points $x_i$ until
they fill up all of $\stackrel{\circ}{\sigma}_i$, we obtain a
deformation retraction $r_1$ from $B_\sigma\backslash
\{x_1,\dots,x_k\}$ to $B_\sigma\backslash
(\stackrel{\circ}{\sigma}_1\cup\dots\cup\stackrel{\circ}{\sigma}_k)$
(see Figure \ref{retractions}).

\item There is another deformation retraction $r_2$ from
$B_\sigma\backslash \{x_1,\dots,x_k\}$ to a subset $S$ of
$\real^d$ that is homeomorphic to the wedge of spheres
$\bigvee^k_{i=1}\mathbb{S}^{d-1}$ (see Figure \ref{retractions}
again). Let $\jmath:\bigvee^k_{i=1}\mathbb{S}^{d-1}\to S$ be a
homeomorphism. Then let $\jmath_i:\sphere^{d-1}\to
S\subset\real^d$ be the maps such that
$\jmath=\jmath_1\vee\dots\vee\jmath_k$.

\item Let $\tilde{\jmath}_i$ be the inclusion
$\sphere^{d-1}\hookrightarrow\bigvee_{i=1}^k\sphere^{d-1}$ of the
$i^\mathrm{th}$ summand of the wedge sum. Furthermore, let
$\mathrm{proj}_i$ be the projection
$$\mathrm{proj}_i:\bigoplus_{i=1}^k\pi_{d-1}(\sphere^{d-1})\twoheadrightarrow\pi_{d-1}(\sphere^{d-1})$$
on the $i^\mathrm{th}$ summand. For $d\geq 3$, the group
homomorphism
$$\tilde{\jmath}:=(\tilde{\jmath}_1)_\ast\circ\mathrm{proj}_1+\dots+(\tilde{\jmath}_k)_\ast\circ\mathrm{proj}_k:\bigoplus_{i=1}^k\pi_{d-1}(\sphere^{d-1})\to\pi_{d-1}(\bigvee_{i=1}^k\sphere^{d-1})$$
is an isomorphism, see Hatcher \cite[Example 4.26]{hat02}. For
$d=2$, this is false. See Remark \ref{d2rem}.

\item Each $(\jmath_i)_\ast$ is an isomorphism, since each
$\jmath_i$ is a homotopy equivalence.

\item This holds because each $I_i$ is an isomorphism.
\end{enumerate}

The diagram is commutative. For every continuous map
$$c:\sphere^{d-1}\to B_\sigma\backslash
(\stackrel{\circ}{\sigma}_1\cup\dots\cup\stackrel{\circ}{\sigma}_k),$$
the equivalence class $[c]$ in $\pi_{d-1}(B_\sigma\backslash
(\stackrel{\circ}{\sigma}_1\cup\dots\cup\stackrel{\circ}{\sigma}_k))$
is mapped by these isomorphisms to
$(W(c,x_1),\dots,W(c,x_k))\in\bigoplus_{i=1}^k\mathbb{Z}$. The
equivalence class $[f|_{\partial\sigma}]$ is therefore mapped to
$(W(f,x_1),\dots,W(f,x_k))=(0,\dots,0)$. Hence
$[f|_{\partial\sigma}]=0$ already in
$\pi_{d-1}(B_\sigma\backslash
(\stackrel{\circ}{\sigma}_1\cup\dots\cup\stackrel{\circ}{\sigma}_k))$.

\begin{figure}
\begin{center}
\psfrag{1}{$r_1$} \psfrag{2}{$r_{2a}$} \psfrag{3}{$r_{2b}$}
\psfrag{x}{$x_i$} \psfrag{B}{$b_\sigma(\sphere^{d-1})$}
\includegraphics[bb = 40 -20 580 830, scale=.65,
clip=true]{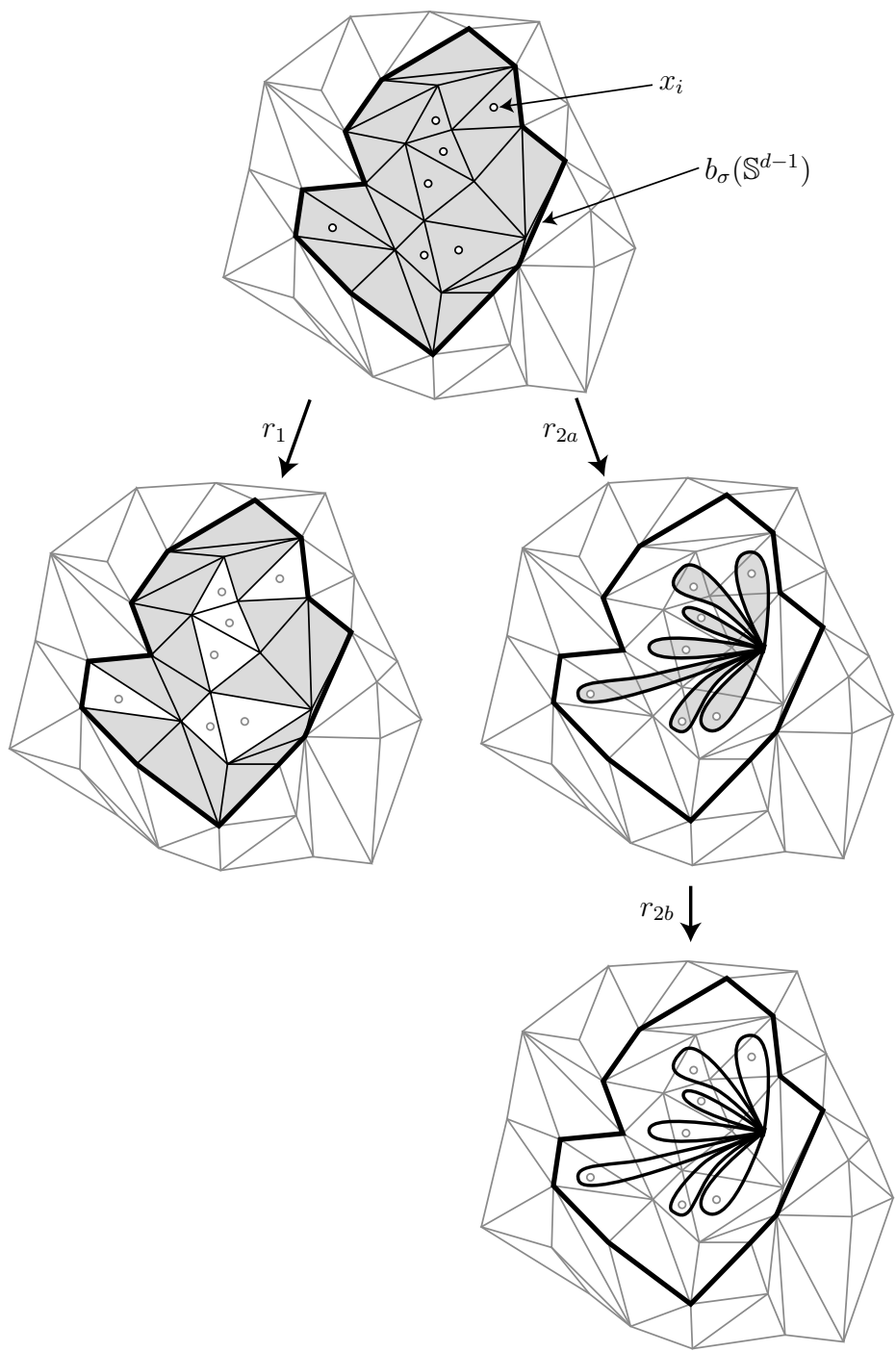}
\end{center}
\caption{The two retractions of
$B_\sigma\backslash\{x_1,\dots,x_k\}$: The dots are the points
$\{x_1,\dots,x_k\}$, the thick line is $b_\sigma(\sphere^{d-1})$
and the thin lines are the $(d-1)$-skeleton of $\Delta_\sigma$.
The left hand side shows the retraction $r_1$ to
$B_\sigma\backslash
(\stackrel{\circ}{\sigma}_1\cup\dots\cup\stackrel{\circ}{\sigma}_k)$,
while the right hand side shows the retraction $r_2$ to
$\bigvee^k_{i=1}\mathbb{S}^{d-1}$, divided into two retractions
$r_{2a}$ and $r_{2b}$.} \label{retractions}
\end{figure}

In other words, $f|_{\partial\sigma}$ is
contractible in $B_\sigma\backslash
(\stackrel{\circ}{\sigma}_1\cup\dots\cup\stackrel{\circ}{\sigma}_k)$.
This means we can extend $f$ continuously to
$\tilde{f}:\df_{d-1}\cup\sigma\to\real^d$ such that
$$\tilde{f}(\sigma)\subset B_\sigma\backslash
(\stackrel{\circ}{\sigma}_1\cup\dots\cup\stackrel{\circ}{\sigma}_k),$$
that is, such that $\tilde{f}(\sigma)$ lies within $\wa$ and the
$(d-1)$-skeleton of the complement of $\wa$. By applying this
argument to all $d$-faces of $\df_d$, we obtain a continuous map
$F:\df_d\to\real^d$.

\emph{Step 3:} We prove that every \tpa\ of $F$ is a \wpa\ of
$f$. Let $P\in\real^d$ be a \tpo\ and
$\sigma_1,\dots,\sigma_q\subset\df_d$ a \tpa\ for $F$. We now
show that $\sigma_1,\dots,\sigma_q$ is also a \wpa\ for $f$ with
\wpo\ $P$:

\begin{itemize}
\item $\dim(\sigma_j)\leq d-1$: In that case we immediately have
$P\in F(\sigma_j)=f(\sigma_j)$.

\item $\dim(\sigma_j)=d$: Suppose $W(f|_{\partial\sigma_j},P)=0$.
For $1\leq i\leq q$, let $\tilde{\sigma}_i$ be the face of
$\Delta_{\sigma_i}$ that contains $P$ in its interior, i.e., the
minimal face containing $P$. We have
\begin{eqnarray*}
d & \stackrel{(1)}{\geq} & \codim(\bigcap^{q}_{i=1}\tilde{\sigma}_i) \\
& \stackrel{(2)}{\geq} & \sum^q_{i=1}\max(0,(d-\dim(\tilde{\sigma}_i))) \\
& = & \sum^q_{i=1}(d-\dim(\tilde{\sigma}_i))\\
& \stackrel{(\ast)}{\geq} & \sum^q_{i=1}(d-\dim(\sigma_i))\\
& = & qd-((d+1)(q-1)+1-q)\\
& = & d.
\end{eqnarray*}

\renewcommand{\labelenumii}{(\theenumii):}
\begin{enumerate}
\item This holds because $\bigcap^{k}_{i=1}\tilde{\sigma}_i$
contains $P$ and therefore is not empty.

\item This holds because the $\Delta_{\sigma_i}$ are in \gp.
\end{enumerate}

The inequality $(\ast)$ is an equality if and only if
$\dim(\tilde{\sigma}_i)=\dim(\sigma_i)$ for all $i$ and in
particular for $i=j$. Hence
$\dim(\tilde{\sigma}_j)=\dim(\sigma_j)=d$. Outside of
$\wa(f|_{\partial\sigma})$, the image $F(\sigma_i)$ lies entirely
in the $(d-1)$-skeleton of $\Delta_{\sigma_i}$; therefore $P$
must lie in $\wa(f|_{\partial\sigma})$.
\end{itemize}
\vskip -.5cm
\end{proof}

\begin{rem}
We used the \dsc\ for continuous maps. This is necessary because
we can not bring a \pl\ approximation of $F$ into \gp\ such that
$F(\sigma)\subset
B_\sigma\backslash(\stackrel{\circ}{\sigma}_1\cup\dots\cup\stackrel{\circ}{\sigma}_k)$.
\end{rem}

\begin{rem}\label{d2rem}
The problem with the case $d=2$ is that $\pi_1$ is not abelian.
Instead of
$$\pi_{d-1}(\real^d\backslash
\{x_1,\dots,x_k\})\cong
\pi_{d-1}(\bigvee^k_{i=1}\mathbb{S}^{d-1})\cong\bigoplus^k_{i=1}\pi_{d-1}(\mathbb{S}^{d-1})\cong\mathbb{Z}^k,$$
which holds for $d \geq 3$, we have
$$\pi_1(\real^2\backslash \{x_1,\dots,x_k\})\cong
\pi_1(\bigvee^k_{i=1}\mathbb{S}^1)\cong F_k \ncong
\bigoplus^k_{i=1}\pi_1(\mathbb{S}^1)$$ where $F_k$ is the free
group on $k$ generators. In particular, $f|_{\partial\sigma}$
need not be contractible in $\real^2\backslash
\{x_1,\dots,x_k\}$, see Figure \ref{d2counterexample}.
\end{rem}

\subsection{The case \emph{d} = 2 of the \wnc}\label{d=2}

We will not show that the cases $d=2$ of the \wnc\ and the \dsc\
are equivalent. Instead, we take a different route:

\begin{prop}\label{wncd2}
If the \wnc\ holds for $d+1$, then it also holds for $d$.
\end{prop}

\begin{cor}
The case $d=3$ of the \dsc\ implies the case $d=2$ of the \wnc.
\end{cor}

\begin{cor}
The \dsc\ implies the \wnc.
\end{cor}

\begin{proof}[Proof of Proposition \ref{wncd2}]
The idea of this proof is based on the proof of Proposition
\ref{tttFromDToD-1} presented in \cite{deL01}.

From any continuous map
$f:\Delta^{(d+1)(q-1)}_{d-1}\rightarrow\real^d$, we construct a
continuous map $F:\Delta^{(d+2)(q-1)}_{d}\rightarrow\real^{d+1}$.
Regard $\real^d$ as the set of all points in $\real^{d+1}$ that
have last coordinate zero. Furthermore, regard
$\Delta^{(d+1)(q-1)}_{d-1}$ as a face of
$\Delta^{(d+2)(q-1)}_{d}$. We denote the vertices of
$\Delta^{(d+1)(q-1)}_{d-1}$ with $e_0,e_1,\dots,e_{(d+1)(q-1)}$
and the $q-1$ additional vertices of $\Delta^{(d+2)(q-1)}_{d}$
with $P_1,P_2,\dots,P_{q-1}$. Now choose $q$ points
$Q,Q_1,Q_2,\dots,Q_{q-1}$ in $\real^{d+1}$, such that $Q$ is
below $\real^d$ (i.e., in $\real^d\times\real^-$) and
$Q_1,\dots,Q_{q-1}$ are above $\real^d$ (i.e., in
$\real^d\times\real^+$). The points $Q_i$ need not be linearly
independent.

Define $F|_{\Delta^{(d+1)(q-1)}_{d-1}}:=f$ and $F(P_i):=Q_i$.
Extend this to all faces of $\Delta^{(d+2)(q-1)}_d$ containing at
least one of the $P_i$ by taking cones over
$F|_{\Delta^{(d+1)(q-1)}_{d-1}}$ with the $Q_i$ as their tips.
More precisely, for all nonnegative numbers $t_i,s_i$ satisfying
$\sum_{i=0}^{(d+1)(q-1)}t_i+\sum_{i=1}^{q-1}s_i=1$, define
$$F\left(\sum^{(d+1)(q-1)}_{i=0}t_ie_i+\sum^{q-1}_{i=1}s_iP_i\right):=
\left(\sum^{(d+1)(q-1)}_{i=0}t_i\right)
f\left(\frac{\sum^{(d+1)(q-1)}_{i=0}t_ie_i}{\sum^{(d+1)(q-1)}_{i=0}t_i}\right)+\sum^{q-1}_{i=1}s_iQ_i.$$
If $\sum_{i=0}^{(d+1)(q-1)}t_i=0$, then the first summand on the
right hand side has to be omitted. Extend this further to all
$d$-faces containing none of the $P_i$ (these are the $d$-faces
of $\df_d$) by using their barycentres to take the cone over
$F|_{\Delta^{(d+1)(q-1)}_{d-1}}$ with $Q$ as its tip (see Figure
\ref{suspension}).

\begin{figure}
\begin{center}
\includegraphics[bb = 80 240 530 655, scale=.6,clip]{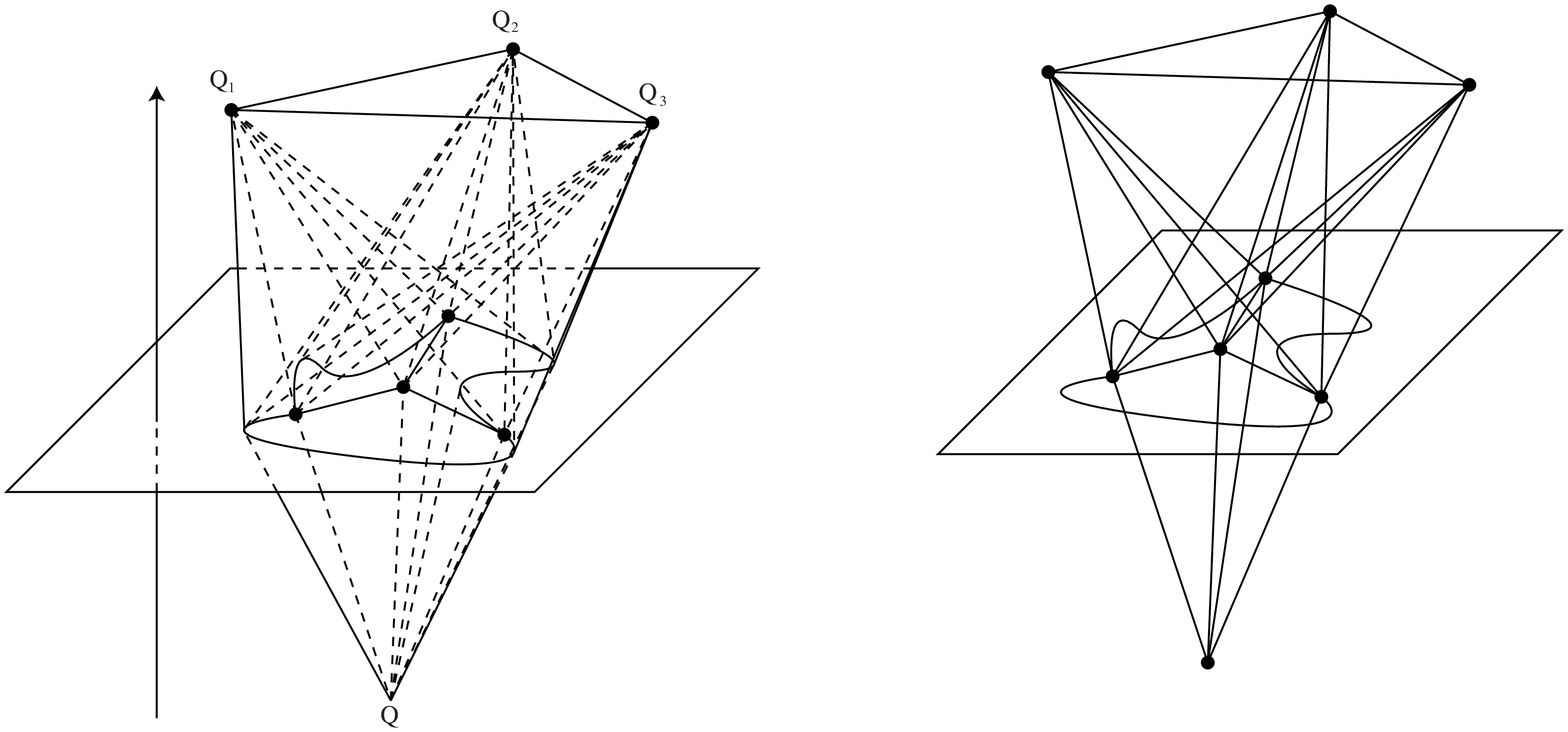}
\end{center}
\caption{The map $F$. The plane $\real^2$ contains the image of
$\Delta^4_1=K_4$, the three points above the plane are
$Q_1,Q_2,Q_3$ and the point below is $Q$.} \label{suspension}
\end{figure}

The \wnc\ for $d+1$ applied to $F$ gives us a \wpo\ $P$ in
$\real^{d+1}$ with a \wpa\ consisting of $q$ disjoint faces
$\sigma_1,\dots,\sigma_q$ of $\Delta^{(d+2)(q-1)}_{d+1}$. If $P$
were above $\real^d$, then all of the $F(\sigma_i)$ would have to
be at least partially above $\real^d$, therefore all of the
$\sigma_i$ would have to contain at least one of the $P_i$. But
this can not be, since the $\sigma_i$ are disjoint, and there are
only $q-1$ points $P_i$. If $P$ were below $\real^d$, then all of
the $F(\sigma_i)$ would have to be at least partially below
$\real^d$, hence all of the $\sigma_i$ would have to contain
$d+1$ of the vertices of $\Delta^{(d+1)(q-1)}_{d+1}$. This cannot
be either, since the $\sigma_i$ are disjoint and there are only
$(d+1)(q-1)+1<(d+1)q$ vertices of $\Delta^{(d+1)(q-1)}_{d+1}$.
Therefore $P$ has to be in $\real^d$. Define
$\tilde{\sigma}_i:=\sigma_i\cap\df_d$. Then
$\tilde{\sigma}_1,\dots,\tilde{\sigma}_q$ are $q$ disjoint faces
that form a \wpa\ for $f$. To see this, we differentiate three
cases.
\begin{itemize}
\item $\dim(\sigma_i)\leq d-1$: $P$ is in $f(\tilde{\sigma}_i)$,
since
$$P\in
F(\sigma_i)\cap\real^d=F(\sigma_i\cap\Delta^{(d+1)(q-1)}_{d-1})=F(\tilde{\sigma}_i)=f(\tilde{\sigma}_i)$$

\item $\dim(\sigma_i)=d$: $P$ is in
$f(\partial\tilde{\sigma}_i)$, since
$$P\in
F(\sigma_i)\cap\real^d=F(\sigma_i\cap\Delta^{(d+1)(q-1)}_{d-1})=F(\partial\tilde{\sigma}_i)=f(\partial\tilde{\sigma}_i)$$

\item $\dim(\sigma_i)=d+1$: W.l.o.g. we may assume that $P$ is
not in $F(\partial\sigma_i)$. We know that $P$ lies in
$\wa(F|_{\partial\sigma_i})\cap\real^d$, therefore
$F(\partial\sigma_i)$ must contain points both above and below
$\real^d$. Thus $\sigma_i$ contains exactly one of the $P_j$,
$\tilde{\sigma}_i$ is $d$-dimensional and we get
\begin{eqnarray*}
P \in \wa(F|_{\partial\sigma_i}) & = & \{x\in\real^{d+1}\mid
W(F|_{\partial\sigma_i},x)\neq
0\}\cap\real^d\\
& =& \{x\in\real^d\mid W(f|_{\partial\tilde{\sigma}_i},x)\neq
0\}\\
& = & \wa(f|_{\partial\tilde{\sigma}_i})
\end{eqnarray*}
\end{itemize}

\end{proof}

\subsection{The connection between \tpa s and \wpa s}
Again, we have proved a stronger statement than just the
equivalence of the \dsc\ and the \wnc.

\begin{thm}\label{ConnectionDSCandWNC}
Let $F:\df_d\to\real^d$ (or even $F:\df\to\real^d$) be a
continuous map. Every \wpa\ of $F|_{\df_{d-1}}$ is a \tpa\ of $F$.

Let $f:\df_{d-1}\to\real^d$ be a continuous map. If $d\geq3$,
then we can extend a slightly distorted version of $f$ to a
continuous map $F:\df_d\to\real^d$ (and even to
$F:\df\to\real^d$) such that every \tpa\ of $F$ is a \wpa\ of
$f$.

If $d=2$, then we can extend a slightly distorted version of $f$
to a continuous map $F:\Delta^{4(q-1)}_3\to\real^3$ (and even to
$F:\Delta^{4(q-1)}\to\real^3$) such that for every \tpa\
$\{\sigma_1,\dots,\sigma_q\}$ of $F$, the set
$\{\sigma_1\cap\Delta^{3(q-1)},\dots,\sigma_q\cap\Delta^{3(q-1)}\}$
is a \wpa\ of $f$.
\end{thm}

\begin{proof}
The first part follows by the proof of Proposition \ref{WNCtoDSC}.

For the second part, we first approximate $f$ by a \pl\ map
$\bar{f}:\df_{d-1}\to\real^d$ in \gp\ without adding any new \wpa
s (possible by Lemma \ref{wncpart2}). We then extend $\bar{f}$ to
the map\linebreak $F:\df_d\to\real^d$ constructed in the proof of
Proposition \ref{wnckgeq3}. From that proof we know that every
\tpa\ of $F$ is a \wpa\ of $\bar{f}$ and thus of $f$.

For the third part, extend $f$ to a continuous map
$\tilde{f}:\Delta^{4(q-1)}_2$ by the suspension described in the
proof of Proposition \ref{wncd2} and proceed with $\tilde{f}$
like we did in the second part.
\end{proof}

This theorem is the correct formulation of the two statements we
speculated about in Remark \ref{BasicIdea}.

\begin{cor}\label{ValidityWNC}
The \wnc\ is valid if $d=1$ and if $q$ is a prime power.
\end{cor}

Therefore Table \ref{ProvenCasesOfTheTTT} applies to the \wnc,
too.

\section{The number of \wpa s and \tpa s}

\begin{prop}[The case \emph{d}$=$1]\label{wncFuerD=1}
For every continuous mapping\linebreak
$f:\Delta^{2(q-1)}_0\to\real$, there are at least $(q-1)!$ \wpa
s. For every continuous map $f:\Delta^{2(q-1)}_1\to\real$
respectively $f:\Delta^{2(q-1)}\to\real$, there are at least
$(q-1)!$ \tpa s.
\end{prop}

\begin{proof}
$\Delta^{2(q-1)}_0$ is a set of $2(q-1)+1=2q-1$ vertices.
$f(\Delta^{2(q-1)}_0)$ is a set of $2(q-1)+1$ real numbers
(counted with multiplicity). Denote the points of
$\Delta^{2(q-1)}_0$, ordered by their function value, by
$P_1,\dots,P_{q-1},M,Y_1,\dots,Y_{q-1}$. A partition of these
points into $q$ sets is a \wpa\ for $f$ if one of the sets is
$\{M\}$ and all the other sets contain exactly one of the $P_i$
and one of the $Q_j$. There are $(q-1)!$ such partitions. The
statement about the number of \tpa s follows directly, because
every \wpa\ of $f|_{\Delta^{2(q-1)}_0}$ is a \tpa\ of $f$.
\end{proof}

\begin{prop}[The case \emph{d}$\geq$3]\label{Equivalence of lower bounds}
If $d\geq3$, then the following three numbers are equal.

\begin{itemize}
\item The minimal number of \tpa s for a continuous map\linebreak
$f:\df\to\real^d$.

\item The minimal number of \tpa s for a continuous map\linebreak
$f:\df_d\to\real^d$.

\item The minimal number of \wpa s for a continuous map\linebreak
$f:\df_{d-1}\to\real^d$.
\end{itemize}

If $d=2$, then at least the first two of these numbers are equal.
\end{prop}

\begin{proof}
Proposition \ref{ConnectionTTTandDSC} and Theorem
\ref{ConnectionDSCandWNC}.
\end{proof}

If Sierksma's conjecture on the minimal number of \tpa s is
correct, then the equivalence established in the previous
proposition carries over to the case $d=2$:

\begin{thm}
The following three statements are equivalent:

\begin{enumerate}
\item Sierksma's conjecture: For all positive integers $d$ and
$q$ and every continuous map $f:\df\to\real^d$ there are at least
$((q-1)!)^d$ \tpa s.

\item For every continuous map $f:\df_d\to\real^d$ there are at
least \linebreak $((q-1)!)^d$ \tpa s.

\item For every continuous map $f:\df_{d-1}\to\real^d$ there are
at least \linebreak $((q-1)!)^d$ \wpa s.
\end{enumerate}
\end{thm}

\begin{proof}
Again by Proposition \ref{ConnectionTTTandDSC} and Theorem
\ref{ConnectionDSCandWNC}, we know that Statements 1 and 2 are
equivalent and that Statement 3 implies Statement 2, which in
turn guarantees Statement 3 if $d\neq 2$.

We now prove that the case $d=3$ of Statement 3 implies the case
\linebreak $d=2$. By Lemma \ref{wncpart2}, it is sufficient to
examine \pl\ maps $f:\Delta^{3(q-1)}_1\to\real^2$ in \gp. Regard
$f$ as a map $\Delta^{3(q-1)}_1\to\real^3$ in the way we did in
the proof of Proposition \ref{wncd2}. For each pair $e_1,e_2$ of
1-dimensional faces of $\Delta^{3(q-1)}_1$, define one of them to
be the ``upper'' and the other one to be the ``lower'' one of the
pair. Now alter $f$ in the following way: For each intersection
$P\in f(e_1)\cap f(e_2)$ of the images of two lines, change $f$
slightly so that the image of the ``upper'' line runs above the
image of the ``lower'' line at $P$, i.e., has a bigger last
coordinate (see Figure \ref{fTofTilde}). We call this new map
$\tilde{f}:\Delta^{3(q-1)}_1\to\real^3$.

\begin{figure*}
\begin{center}
\includegraphics[bb = -5 325 640 515, scale = .6,
clip]{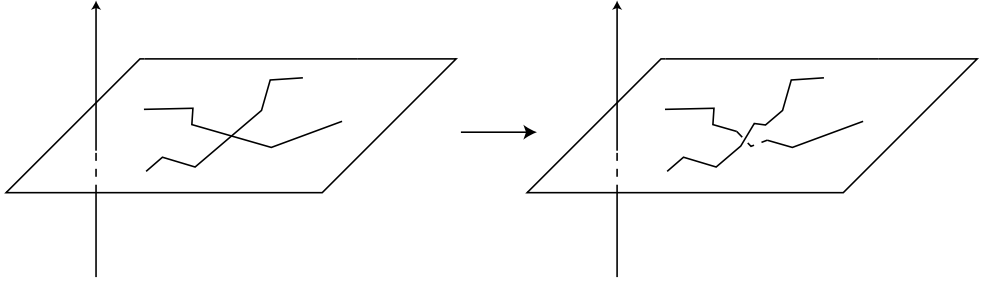}
\end{center}
\caption{How $\tilde{f}$ is obtained from $f$} \label{fTofTilde}
\end{figure*}

We continue similar to the proof of Proposition \ref{wncd2} and
choose points $Q_1,\dots,Q_{q-1}$ high above $\real^2$ and a
point $Q$ far below $\real^2$ and extend $\tilde{f}$ to a map
$F:\Delta^{4(q-1)}_2\to\real^3$ by taking cones using the $Q_i$
and $Q$. Let $\{\sigma_1,\dots,\sigma_q\}$ be a \wpa\ for $F$ and
denote $\tilde{\sigma}_i:=\sigma_i\cap\Delta^{3(q-1)}_1$. By the
argument given in that proof,
$\{\tilde{\sigma}_1,\dots,\tilde{\sigma}_q\}$ is a \wpa\ for $f$.
Since $f$ is in \gp, there are two possibilities for the
$\tilde{\sigma}_i$.
\begin{itemize}
\item All but one of the $\tilde{\sigma}_i$ are 2-dimensional.
The one that is not 2-dimen\-sional, say $\tilde{\sigma}_1$, is
therefore 0-dimensional. Since $\{\sigma_1,\dots,\sigma_q\}$ is a
\wpa\ for our constructed $F$, the faces
$\sigma_2,\dots,\sigma_q$ have to be 3-dimensional and the face
$\sigma_1$ has to be 0-dimensional. Therefore each of the faces
$\sigma_2,\dots\sigma_q$ contains exactly one of the vertices
$P_i$. Hence the \wpa\
$\{\tilde{\sigma}_1,\dots,\tilde{\sigma}_q\}$ of $f$ corresponds
to $(q-1)!$ \wpa s of $F$.

\item All but two of the $\tilde{\sigma}_i$ are 2-dimensional.
The ones that are not \linebreak 2-dimensional, say
$\tilde{\sigma}_1$ and $\tilde{\sigma}_2$, are therefore
1-dimensional. W.l.o.g. let $\tilde{\sigma}_1$ be the ``upper''
one of the two. Since $\{\sigma_1,\dots,\sigma_q\}$ is a \wpa\
for $F$, the faces $\sigma_3,\dots,\sigma_q$ have to be
3-dimensional, the face $\sigma_2$ has to be 2-dimensional and
the face $\sigma_1$ has to be 1-dimensional. Hence the \wpa\
$\{\tilde{\sigma}_1,\dots,\tilde{\sigma}_q\}$ of $f$ corresponds
to $(q-1)!$ \wpa s of $F$.
\end{itemize}

In any case, a \wpa\ of $f$ corresponds to $(q-1)!$ \wpa s of
$F$. Since there are at least $((q-1)!)^3$ \wpa s of $F$, there
have to be at least $((q-1)!)^2$ \wpa s of $f$, which is the
bound we wanted to obtain.
\end{proof}

We now know that the proved and conjectured lower bounds for the
number of \tpa s given in Table \ref{SierksmaNumberTable} also
apply to the number of \wpa s -- except for the proved bound in
the case $d=2$.

\begin{prop}[The case \emph{d}$=$2]
Let $q$ be a prime. There are at least
$$\frac{1}{((q-1)!)^2}\cdot\Big(\frac{q}{2}\Big)^{2(q-1)}$$
\wpa s for every map $f:\Delta^{3(q-1)}_1\to\real^2$.
\end{prop}

\begin{proof}
In the case $d=3$, there are at least
$b:=\frac{1}{(q-1)!}\cdot(\frac{q}{2})^{2(q-1)}$ \tpa s (Theorems
\ref{ManyTverbergPartitions} and \ref{Equivalence of lower
bounds}) and thus the same number of \wpa s. By the proof of the
previous theorem, $\frac{b}{(q-1)!}$ is a bound for the number of
\wpa s for $d=2$.
\end{proof}

The proved and conjectured bounds for $d=2$ are compared in Table
\ref{ComparisonBoundsWNCandTTT}.

\begin{table}
\begin{tabular}{|c||c|c|c|}
  \hline
    & Bound for the number of  & Bound for number of  & Sierksma \\
$q$ & \wpa s & \tpa s & conjecture\\
    & $\frac{1}{((q-1)!)^2}\cdot(\frac{q}{2})^{2(q-1)}$ & $\frac{1}{(q-1)!}\cdot(\frac{q}{2})^{\frac{3}{2}(q-1)}$
 & $((q-1)!)^2$\\
  \hline
  \hline
2  & 1 &  1 &  1\\
  \hline
3  & 2 &  2  & 4\\
  \hline
5  & 3 &  11  &576\\
  \hline
7  & 7 &  110 &518400\\
  \hline
11 & 49 & 35130&   1,31682$\cdot 10^{13}$\\
  \hline
13 & 142 &895579&  2,29443$\cdot 10^{17}$\\
  \hline
17 & 1260 &   967018146   &4,37763$\cdot 10^{26}$\\
  \hline
19 & 3850  &  39101761511 &4,09904$\cdot 10^{31}$\\
  \hline
23 & 37083  & 8,95905$\cdot 10^{13}$ &1,26338$\cdot 10^{42}$\\
  \hline
29 & 1170379 &1,96479$\cdot 10^{19}$ &9,29569$\cdot 10^{58}$\\
  \hline
\end{tabular}
\caption{A comparison of bounds for the number of \tpa s and \wpa
s in the case $d=2$. The numbers are rounded up.}
\label{ComparisonBoundsWNCandTTT}
\end{table}

\begin{ex}
For the alternating linear model of $K_n$ described in Example
\ref{K_nByGuy}, there are $((q-1)!)^2$ \wpa s, exactly the bound
conjectured in the previous Theorem.
\end{ex}

\chapter{\emph{Q}-\w\ Graphs}
The \wnc\ for $d=2$ claims that complete graphs have a certain
property. We now discuss which other graphs have this property,
too.

\begin{defn}
We call a graph $G$ \defw{\emph{q}-\w} if for every map
$f:G\to\real^2$ there are $q$ disjoint paths or cycles
$P_1,\dots,P_q$ in $G$ with
$$\Big(\bigcap_{P_i\mathrm{\ is\ a\ path}}
f(P_i)\Big) \cap\Big(\bigcap_{P_i\mathrm{\ is\ a\ cycle}}
\wa(f|_{P_i})\Big)\neq\emptyset.$$

In accordance with the definition in the previous section, we
call $P_1,\dots,P_q$ a $q$-\wpa\ for $f$.
\end{defn}

The case $d=2$ of the \wnc\ claims that $K_{3q-2}$ is $q$-\w.

\begin{prop}
If $q$ is a prime power, then $K_{3q-2}$ is $q$-winding.
\end{prop}

We now take a closer look at 1-, 2- and 3-\w\ graphs.

\section{1-\w\ graphs}
Every non-empty graph is 1-winding.

\section{2-\w\ graphs and \dy-operations}
Here are two examples of 2-\w\ graphs:

\begin{center}
\begin{tabular}{cc}
\includegraphics*[bb= 60 540 240 720,scale=.75]{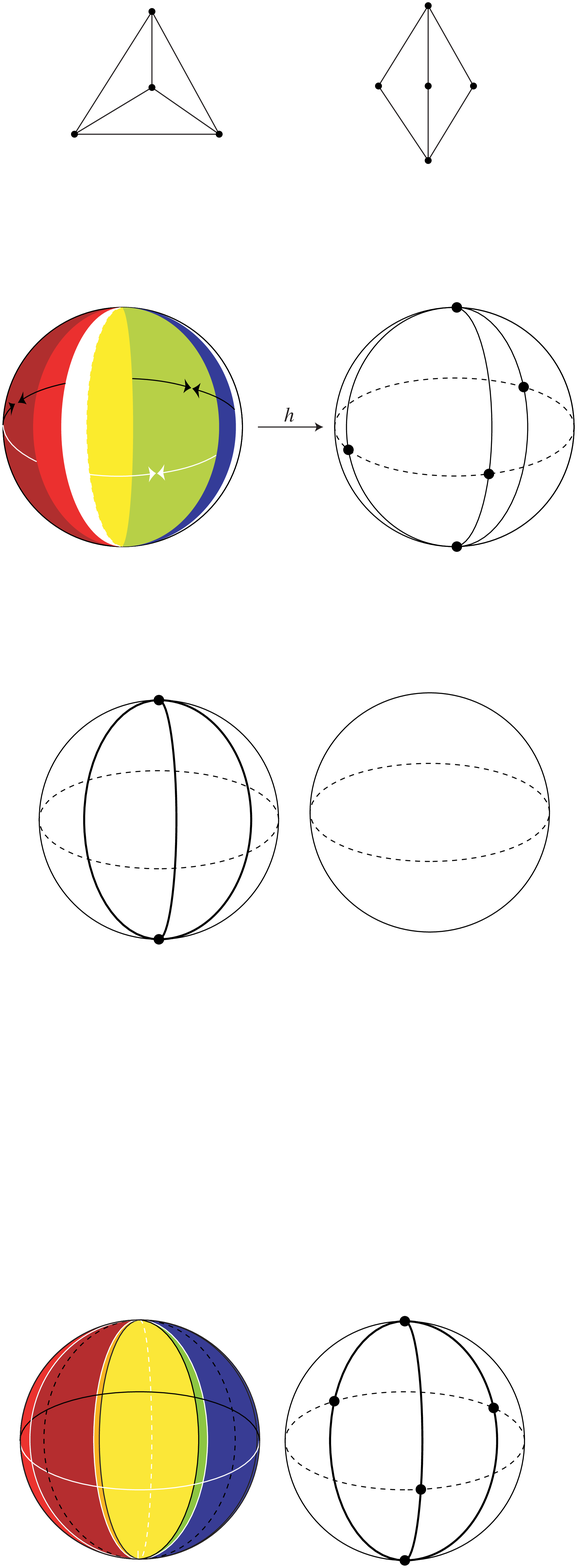}
& \includegraphics*[bb= 370 520 500 710,scale=.75]{K_4UndK_23.eps}  \\
$K_4$ & $K_{2,3}$
\end{tabular}
\end{center}

\begin{prop}\label{K_4ANDK_2,3}
$K_4$ and $K_{2,3}$ are $2$-\w.
\end{prop}

Before we can give a proof, we need to introduce \dy-operations.
We discuss their effect on $q$-winding graphs in general before
we return to the proof of the proposition above.

\begin{defn}
A \dy-operation deletes the three edges of a triangle and adds a
3-valent vertex with edges going from that vertex to the three
vertices of the triangle. A \yd-operation is the reverse of a
\dy-operation.

\begin{figure}[htb]
\begin{center}
\includegraphics*[bb= 4 560 680 830, scale=.55]{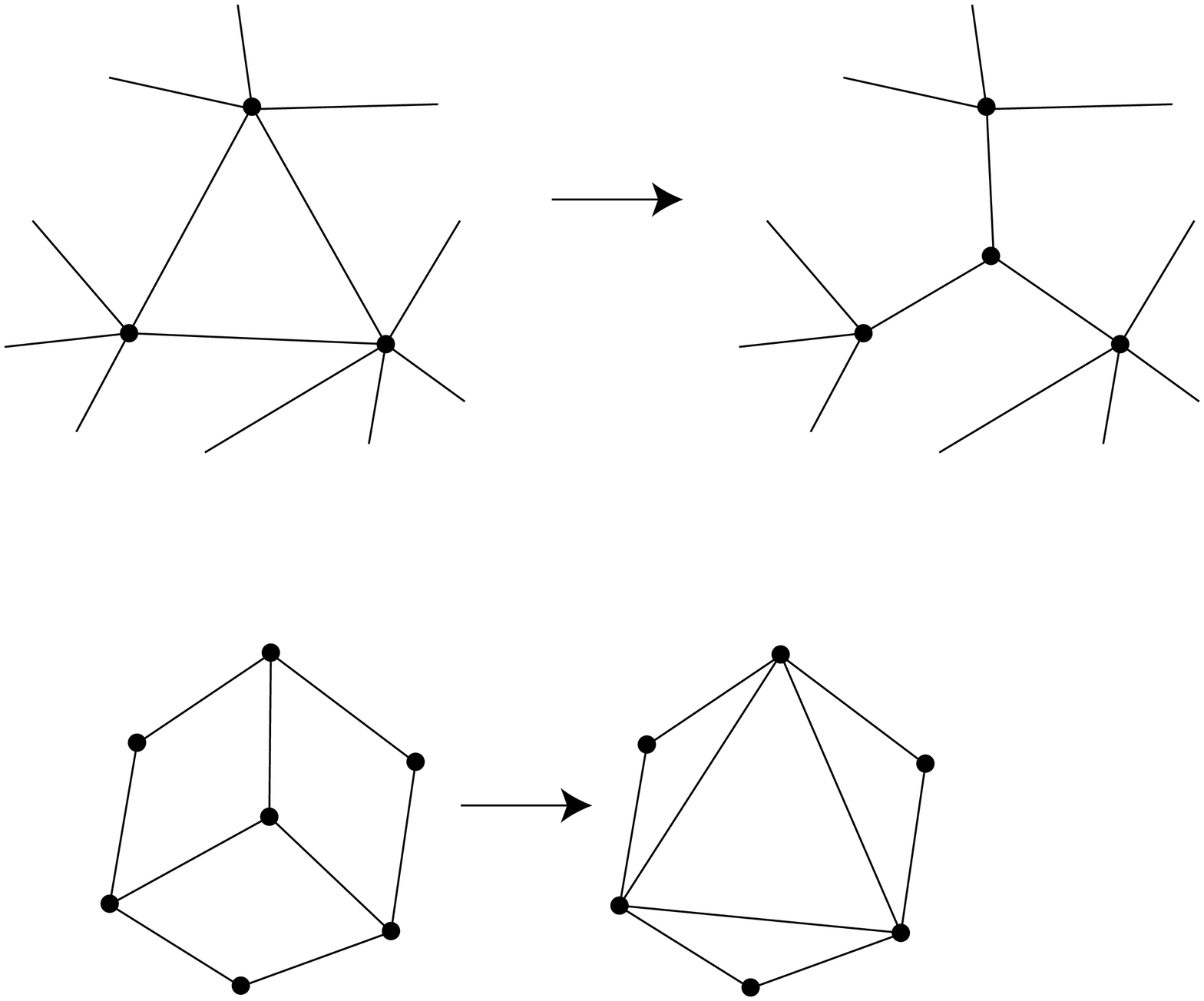}
\end{center}
\caption{An example of a \dy-operation.}
\end{figure}

\end{defn}

\begin{lem}
Let $G$ and $G'$ be two graphs. If there exists a continuous map
$f:G\to G'$ that maps disjoint paths to disjoint paths, then $G'$
is $q$-\w\ if $G$ is $q$-\w.
\end{lem}

\begin{proof}
Let  $g:G'\to\real^2$ be a drawing of $G'$. Then $g\circ
f:G\to\real^2$ is a drawing of $G$. Since $G$ is $q$-\w, there
are $q$ disjoint paths in $G$ that form a $q$-\wpa\ for $g\circ
f$. These paths are mapped under $f$ to $q$ disjoint pathes in
$G'$ that form a $q$-\wpa\ for $g$. Since $g$ was arbitrary, $G'$
is $q$-\w.
\end{proof}

\begin{prop}
A graph obtained from a \dy-operation on a $q$-\w\ graph is again
$q$-\w. A graph obtained from a \yd-operation on a \linebreak
$q$-\w\ graph need not be $q$-\w.
\end{prop}

\begin{figure}[htb]
\begin{center}
\includegraphics*[bb= 4 260 585 470, scale=.7]{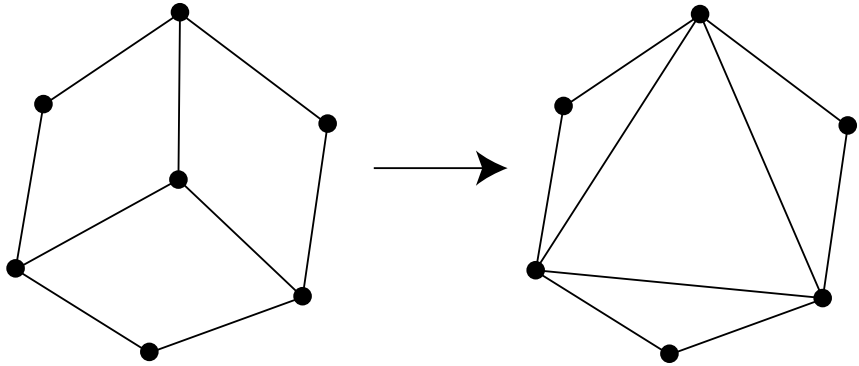}
\end{center}
\caption{A \yd-operation that transforms this 2-\w\ graph into a
graph that is not 2-\w.} \label{YDeltaCounterexample}
\end{figure}

\begin{proof}
Assume that $G'$ is obtained from $G$ by a \dy-operation, more
precisely by deleting the edges $v_1v_2,v_2v_3$ and $v_1v_3$ and
adding the vertex $v$ together with the edges $vv_1,vv_2$ and
$vv_3$. Define $f:G\to G'$ as the identity on all vertices of $G$
and all edges of $G$ except the three deleted ones. For these,
define
\begin{eqnarray*}
f(v_1v_2) & := & v_1vv_2,\\
f(v_2v_3) & := & v_2vv_3,\\
f(v_1v_3) & := & v_1vv_3.
\end{eqnarray*}
The function $f$ maps disjoint paths to disjoint paths. $G'$ is
thus $q$-\w\ if $G$ is $q$-\w.

Figure \ref{YDeltaCounterexample} illustrates that performing
\yd-operations may destroy the property of being $q$-\w.
\end{proof}

We return to the discussion of 2-\w\ graphs.

\begin{proof}[Proof of Proposition \ref{K_4ANDK_2,3}]
The \wnc\ holds for \linebreak $q=2$, hence $K_4$ is 2-\w. The
graph $K_{2,3}$ can be obtained from $K_4$ by a \dy-operation and
hence is 2-\w\ as well.
\end{proof}

\begin{thm}
A graph is $2$-\w\ if and only if it contains $K_4$ or $K_{2,3}$
as a minor.
\end{thm}

\begin{proof}
Every graph that has a $q$-\w\ minor is itself $q$-winding.
Therefore every graph containing $K_4$ or $K_{2,3}$ as a minor is
2-\w.

On the other hand, if a graph does not contain one of these two
graphs as a minor, then it is \defw{outerplanar}, that is, it has
a planar drawing with all vertices lying on the exterior region
(Chartrand and Harary \cite{chh67}). In such a drawing no two
edges intersect (the drawing is planar!) and no cycle winds around
a vertex. Hence the graph is not 2-\w.
\end{proof}

\section{3-\w\ graphs and \emph{q}-\w\ subgraphs of complete graphs}
We prove two general results about $q$-\w\ subgraphs of
$K_{3q-2}$ and obtain the minimal 3-\w\ subgraph of $K_7$.

\begin{thm}\label{wg1}
Let $p\geq 3$ be a prime and $M$ a maximal matching in
$K_{3p-2}$. Then $K_{3p-2}-M$ is $p$-\w.
\end{thm}

\begin{proof}
Let $N:=4(p-1)$ and let $f:K_{3p-2}\to\real^2$ be a drawing of
$K_{3p-2}$. We divide the proof in three steps.
\begin{enumerate}
\item We describe a $\mathbb{Z}_p$-invariant subcomplex $L$ of
$(\Delta^N)^{*p}_{\Delta(2)}$.

\item We show that $\ind_{\mathbb{Z}_p}(L)\geq
N>N-1=\ind_{\mathbb{Z}_p}((\real^3)^{*p}_{\Delta})$. By Lemma~
\ref{IndexTheoryLemma1} on index theory, $L$ cannot be mapped to
$(\real^3)^{*p}_{\Delta}$ $\mathbb{Z}_p$-equivariantly.

\item We extend the drawing $f$ to a map $F:\Delta^N\to\real^3$
and examine the \tpa s of $F$ and \wpa s of $f$ that correspond
to $L$.
\end{enumerate}

\emph{Step 1:} The maximal simplices of
$(\Delta^N)^{*p}_{\Delta(2)}$ correspond to the edges of a
complete $(N+1)$-partite hypergraph with $p$ vertices in each
shore. In Figure~\ref{dots}, the $N+1$ rows represent the $N+1$
shores. We extend the matching $M$ of $K_{3p-2}$ to a maximal
matching the vertices of $\Delta^N$ and group the rows into pairs
accordingly. One row remains single. For each pair of rows we now
choose a $\mathbb{Z}_p$-invariant cycle in the complete bipartite
graph generated by these two shores, such that the cycles contain
no vertical lines. The maximal simplices of $L$ shall be the
maximal simplices of $(\Delta^N)^{*p}_{\Delta(2)}$ whose
corresponding edge in the hypergraph contains an edge of each
cycle (see Figure \ref{dots}). Through this, $L$ is completely
determined and $\mathbb{Z}_p$-invariant in
$(\Delta^N)^{*p}_{\Delta(2)}$.

\begin{figure}
\begin{center}
\includegraphics[bb = 50 380 540 790, scale=.7, clip]{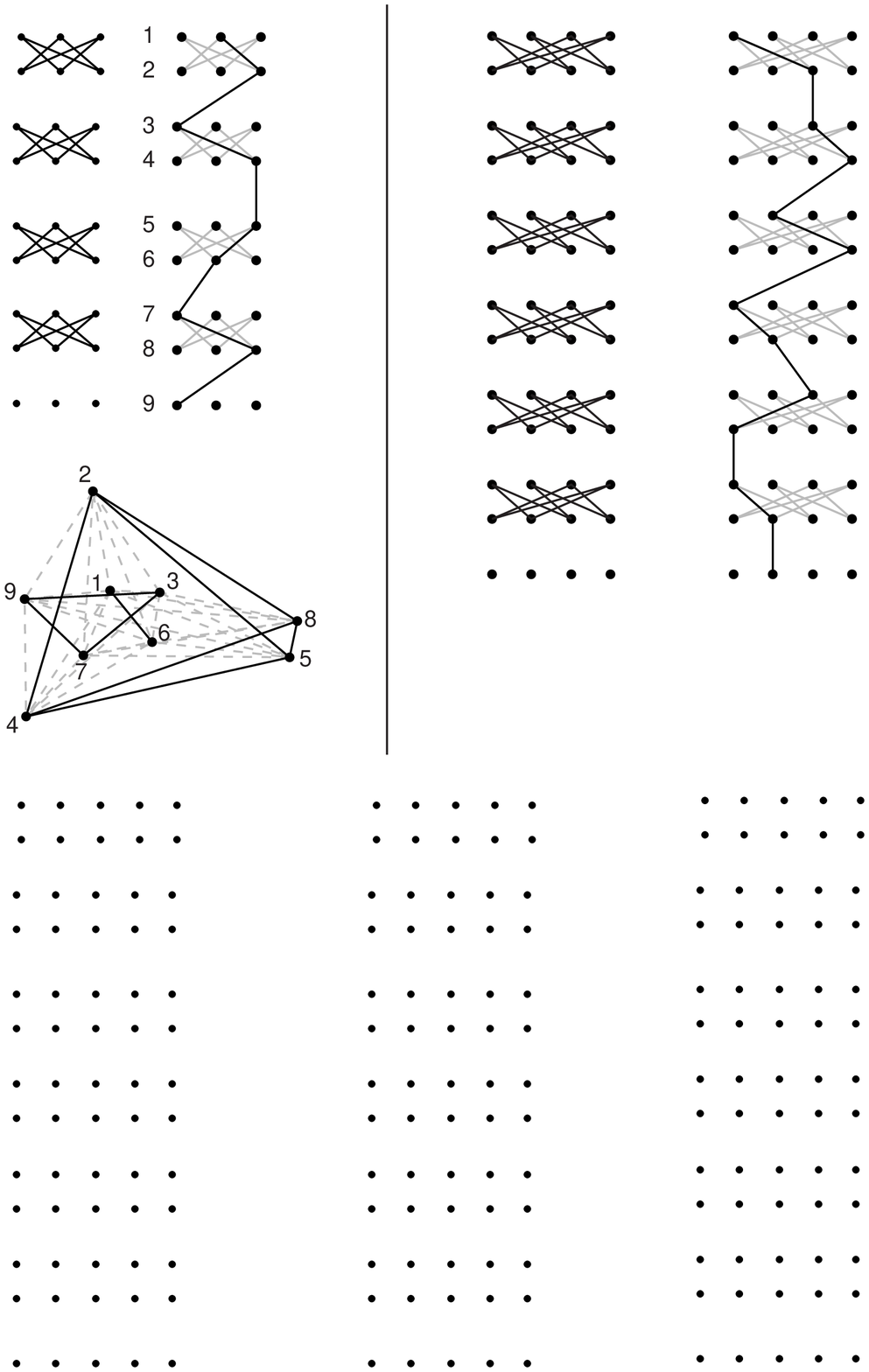}
\end{center}
\caption{This figure illustrates the correspondence between
$(\Delta^N)^{*p}_{\Delta(2)}$ and the complete $(N+1)$-partite
hypergraph with $p$ vertices in each shore. The left hand side
shows the case $p=3$ and $N=8$: The rows represent the 9 shores of
3 vertices each. For each pair of rows a cycle is drawn. The thick
line corresponds to a maximal face of $L$, the partition of the
vertices of $(\Delta^8)^{*3}_{\Delta(2)}$ represented by this face
is drawn in black below the dots. The right hand side of the
figure shows the hypergraph corresponding to
$(\Delta^{13})^{*4}_{\Delta(2)}$.}
\label{dots}
\end{figure}

\emph{Step 2:} $L$ can be interpreted as the join of its $N/2$
circles and the remaining row of $p$ points:
$$L\cong(\sphere^1)^{*N/2}*D_p.$$
By Lemma \ref{IndexTheoryLemma2}, we obtain
\begin{eqnarray*}
\ind_{\mathbb{Z}_p}(L) & = &
\ind_{\mathbb{Z}_p}((\sphere^1)^{*N/2}*D_p)\\
& \geq & \frac{N}{2}\ind_{\mathbb{Z}_p}(\sphere^1)+\frac{N}{2}+\ind_{\mathbb{Z}_p}(D_p)\\
& = & \frac{N}{2}+\frac{N}{2}+0\\
& = & N.
\end{eqnarray*}
The identity $N-1=\ind_{\mathbb{Z}_p}((\real^3)^{*p}_{\Delta})$
was also stated in Lemma \ref{IndexTheoryLemma2}.

\emph{Step 3:} By Theorem \ref{ConnectionDSCandWNC}, we can
extend a slightly distorted version of $f$ to a continuous map
$F:\Delta^{4(q-1)}\to\real^3$, such that for every \tpa\
$\{\sigma_1,\dots,\sigma_q\}$ of $F$, the set
$\{\sigma_1\cap\Delta^{3(q-1)},\dots,\sigma_q\cap\Delta^{3(q-1)}\}$
is a \wpa\ for $f$.

The maximal simplices of $L$ correspond to sets of disjoint faces
of $\Delta^N$. For every continuous map $\Delta^N\to\real^3$, at
least one of these sets is a \tpa, because $L$ can not be mapped
$\mathbb{Z}_p$-equivariantly to $(\Delta^N)^{*q}_{\Delta(2)}$.

Since we chose the cycles in Step 1 such that they contain no
vertical lines, we can be sure that in every such \tpa\
$\{\sigma_1,\dots,\sigma_q\}$ of $F$ the two vertices that form a
pair do not belong to the same face. This also holds for the
corresponding \wpa\
$\{\sigma_1\cap\Delta^{3(q-1)},\dots,\sigma_q\cap\Delta^{3(q-1)}\}$
of $f$. Therefore we can delete the edges in $K_{3p-2}$ connecting
the two vertices of a pair (that is, we can delete the maximal
matching $M$) and still have a $p$-\w\ graph.
\end{proof}

\begin{rem}
The complex $L$ was used before to obtain a lower bound for the
number of Tverberg partitions (Theorem
\ref{ManyTverbergPartitions}, see Matou{\v{s}}ek
\cite[Theorem~6.6.1] {mat03}).
\end{rem}

\begin{prop}\label{wg2}
Let $N$ be $q-1$ edges of $K_{3q-2}$ meeting in one vertex. Then
$K_{3q-2}-N$ is not $q$-\w.
\end{prop}

\begin{proof}
All we need to do is to find a drawing of $K_{3q-2}-N$ without a
\linebreak $q$-\wpa. We can use the alternating linear model of
$K_n$ described in Example \ref{K_nByGuy}. All we have to do is
to order the vertices such that the meeting vertex is at the
right end of the drawing and the other vertices of $N$ have the
numbers $1,3,5,\dots,2q-5,2q-3$. The edges of $N$ are then in the
upper half. Figure \ref{K_nOhneKante} illustrates the situation
for $K_7$ and $K_{10}$.

\begin{figure}
\begin{center}
\includegraphics*[bb= 110 550 470 670,scale=.7]{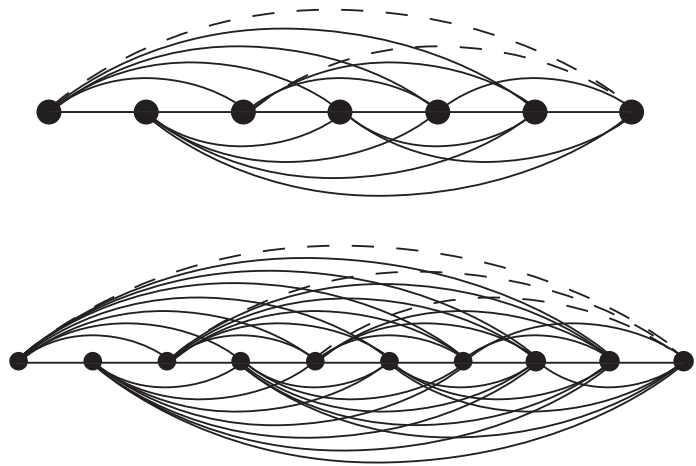}
\includegraphics*[bb= 100 390 500 530,scale=.7]{K_nAlsLinieOhneZweiKanten.eps}
\end{center}
\caption{Drawing of $K_7$ and $K_{10}$. The edges that form $N$
are dashed.} \label{K_nOhneKante}
\end{figure}
\end{proof}

\begin{cor}
The unique minimal $3$-\w\ minor of $K_7$ is $K_7-M$, where $M$
is a maximal matching.
\end{cor}

\begin{proof}
$K_7-M$ is a 3-\w\ minor of $K_7$ (Theorem \ref{wg1}). It is
minimal, because all edges not in $M$ are adjacent to an edge in
$M$ and thus must not be deleted (Proposition \ref{wg2}).

If on the other hand $K$ is 3-\w\ minor of $K_7$, then only a
matching can be deleted (again by Proposition \ref{wg2}). For $K$
to be minimal, this matching must be maximal.
\end{proof}

\begin{prop}
Not every $3$-\w\ graph has $K_7$ minus a maximal matching as a
minor.
\end{prop}

\begin{proof}
Let $M$ be a maximal matching in $K_7$. Execute a \dy-operation
on an $K_7-M$; the resulting graph is 3-\w, but does not have
$K_7-M$ as a minor.
\end{proof}

\clearpage \addcontentsline{toc}{chapter}{Bibliography}
\bibliographystyle{alpha}
\bibliography{Diplom}







\end{document}